# CLASSICAL AND FREE INFINITELY DIVISIBLE DISTRIBUTIONS AND RANDOM MATRICES

By Florent Benaych-Georges

*École Normale Supérieure*

We construct a random matrix model for the bijection $\Psi$ between clas- sical and free infinitely divisible distributions: for every $d \geq 1$, we associate in a quite natural way to each $*$-infinitely divisible distribution $\mu$ a distribution $\mathbb{P}_d^\mu$ on the space of $d \times d$ Hermitian matrices such that $\mathbb{P}_d^\mu * \mathbb{P}_d^\nu = \mathbb{P}_d^{\mu * \nu}$. The spectral distribution of a random matrix with distribution $\mathbb{P}_d^\mu$ converges in probability to $\Psi(\mu)$ when $d$ tends to $+\infty$. It gives, among other things, a new proof of the almost sure convergence of the spectral distribution of a matrix of the GUE and a projection model for the Marchenko–Pastur distribution. In an analogous way, for every $d \geq 1$, we associate to each $*$-infinitely divisible distribution $\mu$, a distribution $\mathbb{L}_d^\mu$ on the space of complex (non-Hermitian) $d \times d$ random matrices. If $\mu$ is symmetric, the symmetrization of the spectral distribution of $|M_d|$, when $M_d$ is $\mathbb{L}_d^\mu$-distributed, converges in probability to $\Psi(\mu)$.

**Introduction.** Free convolution $\boxplus$, defined in Bercovici and Voiculescu (1993), is a binary operation on the set of probability measures on the real line, arising from free probability theory ($\mu \boxplus \nu$ is the distribution of $X + Y$ when $X, Y$ are free and have distributions $\mu, \nu$). It is associative, commutative and continuous with respect to the weak convergence. A probability measure $\mu$ on $\mathbb{R}$ is said to be $\boxplus$-infinitely divisible if for every $n \geq 1$, there exists a probability measure $\mu_n$ on $\mathbb{R}$ such that $\mu_n^{\boxplus n}$ equals to $\mu$.

It is shown in Bercovici, Pata and Biane (1999) that there exists an home-omorphism $\Psi$ from the set of $*$-infinitely divisible distributions to the set of $\boxplus$-infinitely divisible distributions which associates to every classical (resp. free) limit theorem a free (resp. classical) analogue. Indeed, for every $*$-infinitely divisible distribution $\mu$, for every sequence $(\mu_n)$ of probability measures, for every sequence $(k_n)$ of integers tending to infinity, the sequence $\mu_n^{* k_n}$ tends to $\mu$ if and only if the sequence $\mu_n^{\boxplus k_n}$ tends to $\Psi(\mu)$.

The proofs in Bercovici, Pata and Biane ([1999](#)) rely on integral transformations and complex analysis. We will, in this article, construct a matricial model for the $\boxplus$-infinitely divisible distributions, and present in a more palpable way the bijection $\Psi$.

Let $\mu$ be an $*$-infinitely divisible distribution. Let $(\mu_n)$ be a sequence of probability measures and $(k_n)$ a sequence of integers which tends to infinity such that the sequence $\mu_n^{*k_n}$ tends weakly to $\mu$. Let, for $d \geq 1$ and $n \geq 1$, $\mathbb{Q}_d^{\mu_n}$ (resp. $\mathbb{K}_d^{\mu_n}$) be the distribution of $U \operatorname{diag}(X_{n,1}, \ldots, X_{n,d}) U^*$ [resp. of $U \operatorname{diag}(X_{n,1}, \ldots, X_{n,d}) V$], where $U, V$ are independent unitary Haar distributed random matrices independent of the i.i.d. random variables $X_{n,1}, \ldots, X_{n,d}$ with distribution $\mu_n$. We will prove, in Section [3](#) (resp. Section [7.1](#)), that the sequence $((\mathbb{Q}_d^{\mu_n})^{*k_n})$ [resp. $((\mathbb{K}_d^{\mu_n})^{*k_n})$] converges weakly to a probability measure $\mathbb{P}_d^{\mu}$ (resp. $\mathbb{L}_d^{\mu}$). The main results of this article are the following ones: the spectral distribution of a random matrix with distribution $\mathbb{P}_d^{\mu}$ converges in probability to $\Psi(\mu)$ when $d$ tends to infinity, and so does the symmetrization of the spectral distribution of $|M_d|$ when $M_d$ is distributed according to $\mathbb{L}_d^{\mu}$. So we have constructed matrix models which go from $*$-infinitely divisible distributions to $\boxplus$-infinitely divisible distributions when the dimension goes from one to infinity. What is more, for all $*$-infinitely divisible distributions $\mu, \nu$ and all $d$, $\mathbb{P}_d^{\mu*\nu} = \mathbb{P}_d^{\mu} * \mathbb{P}_d^{\nu}$ and $\mathbb{L}_d^{\mu*\nu} = \mathbb{L}_d^{\mu} * \mathbb{L}_d^{\nu}$. This property (and the fact that all formulas depend analytically on $d$, so could be extended to noninteger $d$) opens the perspective of a continuum between the classical convolution $*$ and the free convolution $\boxplus$ for infinitely divisible mesures [M. Anshelevich has already constructed such a continuum in Anshelevich ([2001](#)), but the model we present here does not interpolate his construction]. T. Cabanal-Duvillard, in Cabanal-Duvillard ([2004](#)), has studied at the same time as the author the distributions $\mathbb{P}_d^{\mu}$, and has proved the same result, but with different methods (processes, measure concentration, integral transforms).

At last, in the case where $\mu$ is the standard normal distribution, $\Psi(\mu)$ is the semi-circle distribution with center zero and radius two, and the distribution $\mathbb{P}_d^{\mu}$ is closely related to the one of the GUE, so that the convergence of the spectral distribution of a matrix with distribution $\mathbb{P}_d^{\mu}$ implies Wigner's result. Likewise, the distribution $\mathbb{L}_d^{\mu}$ is the one of a matrix with independent Gaussian entries, and we have a new proof of the convergence of the spectral distribution of the Wishart matrix with parameter 1 to the Marchenko–Pastur distribution.

In the same way, in the case where $\mu$ is the classical Poisson distribution, this result allows us to see the Marchenko–Pastur distribution as the limit spectral distribution of a sum of independent rank-one projections.

The text is organized as follows. In Section [1](#) we recall a few results about infinitely divisible distributions and about their classical and free cumulants. In Section [2](#) we explain the choice of the model (i.e., of the distributions $\mathbb{P}_d^\mu$ and $\mathbb{L}_d^\mu$). In Section 3 we construct the distributions $\mathbb{P}_d^\mu$. Finally, the convergence in probability of the spectral distribution of a random matrix with distribution $\mathbb{P}_d^\mu$ to $\Psi(\mu)$ is proved in two steps. In the first one, we show the convergence when the Lévy measure has compact support, and in the second one (in Section 6), we extend this result using approximation and compound Poisson distributions. The first step is achieved with the moment method, and is divided into two steps: convergence of the mean of every moment in Section 4, almost sure convergence in Section 5. The distributions $\mathbb{L}_d^\mu$ are constructed in Section 7.1, the convergence in probability of the symmetrization of the spectral distribution of $|M_d|$, when $M_d$ is distributed according to $\mathbb{L}_d^\mu$, is also divided in two steps.

## 1. Preliminary results about infinitely divisible distributions.

1.1. *Definitions and the bijection* $\Psi$. The results of this section concerning classical probabilities are in Gnedenko and Kolmogorov (1954) and in Petrov (1995); the results concerning free probabilities are in Bercovici and Voiculescu (1993) and in Bercovici, Pata and Biane (1999), except the continuity of the inverse of the bijection $\Psi$, which is shown in Barndorff-Nielsen and Thorbjørnsen (2002). A probability measure $\mu$ on $\mathbb{R}$ is said to be $*$-infinitely divisible (resp. $\boxplus$-infinitely divisible) if for every $n \geq 1$, there exists a probability measure $\mu_n$ on $\mathbb{R}$ such that $\mu_n^{*n}$ (resp. $\mu_n^{\boxplus n}$) equals $\mu$, which is equivalent to the existence of a sequence $(\mu_n)$ of probability measures, of a sequence $(k_n)$ of integers tending to infinity, such that $\mu_n^{*k_n}$ (resp. $\mu_n^{\boxplus k_n}$) tends weakly to $\mu$.

We can characterize $*$-infinitely divisible distributions (resp. $\boxplus$-infinitely divisible distributions) with their Fourier transform (resp. their Voiculescu transform). A probability measure $\mu$ on $\mathbb{R}$ is $*$-infinitely divisible (resp. $\boxplus$-infinitely divisible) if and only if there exists a real $\gamma$ and a positive finite measure $G$ on $\mathbb{R}$ such that its Fourier transform $\widehat{\mu}$ (resp. its Voiculescu transform $\varphi_\mu$) has the form

$$\widehat{\mu}(t) = \exp\left\{ i\gamma t + \int_{u \in \mathbb{R}} \underbrace{\left[ e^{itu} - 1 - \frac{itu}{1+u^2} \frac{1+u^2}{u^2} \right]}_{=-\frac{t^2}{2} \text{ for } u=0} \mathrm{d}G(u) \right\}$$

$$\left( \text{resp. } \varphi_\mu(z) = \gamma + \int_{t \in \mathbb{R}} \frac{1+tz}{z-t} \mathrm{d}G(t) \right).$$

In this case, the pair $(\gamma, G)$ is unique, and we denote $\mu = \nu_*^{\gamma, G}$ (resp. $\nu_\boxplus^{\gamma, G}$).

REMARK. There exists other parametrizations of $*$-infinitely divisible distributions: for example, denoting $\gamma' = \gamma, \sigma^2 = G(\{0\}), L(A) = \int_A \frac{1+u^2}{u^2} \mathrm{d}G(u)$ for all Borel set $A$ of $\mathbb{R} \setminus \{0\}$, one has $\int_{u \in \mathbb{R} \setminus \{0\}} (1 \wedge u^2) \mathrm{d}L(u) < \infty$, and $\widehat{\mu}(t) = \exp(i\gamma' t - \frac{\sigma^2 t^2}{2} + \int_{u \in \mathbb{R} \setminus \{0\}} (e^{itu} - 1 - \frac{itu}{1+u^2}) \mathrm{d}L(u)).$

We now give the definition of $\Psi$, referred to in the Introduction.

**Theorem 1.1** (Bercovici–Pata's bijection). *We endow the set of positive finite measures on $\mathbb{R}$ with the weak topology; the subsets $\{*$-infinitely divisible distributions$\}$ and $\{\boxplus$-infinitely divisible distributions$\}$ are also endowed with the weak topology.*

1. *The maps*

$$\mathbb{R} \times \{\text{ positive finite measures}\} \to \{*\text{-infinitely divisible distributions}\},$$

$$(\gamma, G) \mapsto \nu_*^{\gamma, G}$$

   *and*

$$\mathbb{R} \times \{\text{ positive finite measures}\} \to \{*\text{-infinitely divisible distributions}\},$$

$$(\gamma, G) \mapsto \nu_{\boxplus}^{\gamma, G}$$

   *are homeomorphisms and we have*

$$\nu_*^{\gamma+\gamma', G+G'} = \nu_*^{\gamma, G} * \nu_*^{\gamma', G'},$$

$$\nu_{\boxplus}^{\gamma+\gamma', G+G'} = \nu_{\boxplus}^{\gamma, G} \boxplus \nu_{\boxplus}^{\gamma', G'}.$$

2. *Let us define the map $\Psi$, from the set of $*$-infinitely divisible distributions to the set of $\boxplus$-infinitely divisible distributions, which maps, for all $(\gamma, G)$, the measure $\nu_*^{\gamma, G}$ to the measure $\nu_{\boxplus}^{\gamma, G}$. Then*

   (a) *$\Psi$ is an homeomorphism called Bercovici–Pata's bijection,*

   (b) *for all $\mu, \nu$ $*$-infinitely divisible distributions, $\Psi(\mu * \nu) = \Psi(\mu) \boxplus \Psi(\nu)$,*

   (c) *Dirac measures are invariant under $\Psi : \Psi(\delta_a) = \delta_a$,*

   (d) *$\Psi(N(m, r^2))$ is the semi-circle distribution with mean $m$ and variance $r^2$, which is $w_{m, 2r}(x)\, \mathrm{d}x$, with*

   $$w_{m, 2r}(x) = \frac{1}{2\pi r^2}(4r^2 - (x - m)^2)^{1/2} \mathbb{1}_{|x-m| \le 2r},$$

   (e) *$\Psi$, restricted to the Cauchy type, is the identity: for all $a > 0$, $\Psi(C_a) = C_a$, where $C_a = \frac{1}{\pi} \frac{a\, \mathrm{d}x}{a^2 + x^2}$,*

   (f) *for all sequence $(\mu_n)$ of probability measures on $\mathbb{R}$, for all sequence $(k_n)$ of integers tending to infinity, the sequence $\mu_n^{*k_n}$ converges weakly to a $*$-infinitely divisible distribution $\mu$ if and only if $\mu_n^{\boxplus k_n}$ converges weakly to $\Psi(\mu)$.*

**Remark 1.2.** In the text, the positive finite measure $G$ is called the Lévy measure of $\nu_*^{\gamma, G}$ and $\nu_{\boxplus}^{\gamma, G}$. We will use the two following properties:

1. If the Lévy measure of a $\boxplus$-infinitely divisible distribution $\nu$ has compact support, then so does $\nu$ [see Hiai and Petz ([2000](#))].

2. $\nu_{\boxplus}^{\gamma, G}$ is symmetric if and only if $\nu_*^{\gamma, G}$ is symmetric, if and only if $G$ is symmetric and $\gamma = 0$.

1.2. *Classical compound Poisson distributions, approximation of ∗-infinitely divisible distributions by ∗-infinitely divisible distributions with compactly supported Lévy measures.*

DEFINITION 1.3. Let $\lambda$ be a nonnegative real, $\rho$ be a probability measure on $\mathbb{R}$. Then the sequence of probability measures on $\mathbb{R}$

$$\left( \left( 1 - \frac{\lambda}{n} \right) \delta_0 + \frac{\lambda}{n} \rho \right)^{*n}, \qquad n \geq 1,$$

converges weakly to a distribution noted $\pi^*_{\rho,\lambda}$, with Fourier transform

$$\widehat{\pi^*_{\rho,\lambda}}(t) = \exp(\lambda(\widehat{\rho}(t) - 1)),$$

where $\widehat{\rho}$ is the Fourier transform of $\rho$.

REMARK 1.4. $\pi^*_{\rho,\lambda}$ is $\nu^{\gamma,G}_*$, with

$$G = \lambda \frac{u^2}{1 + u^2} \, \mathrm{d}\rho(u), \qquad \gamma = \lambda \int_{u \in \mathbb{R}} \frac{u}{1 + u^2} \, \mathrm{d}\rho(u).$$

We introduce now the compactly supported approximations of the positive finite measure $G$.

DEFINITION 1.5. Let, for $G$ positive finite measure on $\mathbb{R}$, $t > 0$, $G^0_t, G_t$ be the positive finite measures on $\mathbb{R}$ defined by

$$G^0_t(A) = G(A \cap [-t, t]), \qquad G_t(A) = G(A \setminus [-t, t])$$

for all Borel set $A$ of $\mathbb{R}$.

We define $\lambda_t \geq 0$, the probability measure $\rho_t$ on $\mathbb{R}$, and $a_t \in \mathbb{R}$ with

$$\lambda_t = \int_{u \in \mathbb{R} \setminus [-t,t]} \frac{1 + u^2}{u^2} \, \mathrm{d}G(u), \qquad \rho_t = \frac{1}{\lambda_t} \frac{1 + u^2}{u^2} \, \mathrm{d}G_t(u),$$

$$a_t = - \int_{u \in \mathbb{R} \setminus [-t,t]} (1/u) \, \mathrm{d}G(u).$$

We will use the following approximation:

$$(2) \qquad \forall \, t > 0 \qquad \nu^{\gamma,G}_* = \nu^{\gamma + a_t, G^0_t}_* * \pi^*_{\rho_t, \lambda_t},$$

because one observes that $\pi^*_{\rho_t, \lambda_t} = \nu^{\alpha, H}_*$ with

$$H = \lambda_t \frac{u^2}{1 + u^2} \, \mathrm{d}\rho_t(u) = G_t, \qquad \alpha = \lambda_t \int_{u \in \mathbb{R}} \frac{u}{1 + u^2} \, \mathrm{d}\rho_t(u) = -a_t.$$

1.3. *Partitions, moments and cumulants of infinitely divisible distributions.* For every probability measure $\mu$, we will denote, when it is defined, by $m_n(\mu)$ the $n$th moment of $\mu$, which is $\int x^n \, d\mu(x)$. In this case, we will denote by $\mathfrak{C}_n(\mu)$ [resp. $\mathfrak{K}_n(\mu)$] its $n$th classical (resp. free) cumulant. Recall that [see Section 4 of Speicher ([1994](#)) or Section 2.5 of Hiai and Petz ([2000](#))]

$$(3) \qquad m_k(\mu) = \sum_{\pi \in \mathrm{Part}(k)} \underbrace{\prod_{V \in \pi} \mathfrak{C}_{|V|}(\mu)}_{\text{denoted by } \mathfrak{C}_\pi(\mu)},$$

$$(4) \qquad m_k(\mu) = \sum_{\pi \in \mathrm{NC}(k)} \underbrace{\prod_{V \in \pi} \mathfrak{K}_{|V|}(\mu)}_{\text{denoted by } \mathfrak{K}_\pi(\mu)},$$

where $\mathrm{Part}(k)$ denotes the set of the partitions of $\{1, \ldots, k\}$ and $\mathrm{NC}(k)$ denotes the set of noncrossing partitions of $[k] = \{1, \ldots, k\}$ (a *noncrossing partition* of a finite totally ordered set $I$ is a partition $\pi$ of $I$ such that there does not exist $x < y < z < t \in I$ with $x$ and $z$ belonging to the same class and $y$ and $t$ belonging to another class).

We will need the following proposition [part of which was proved in Barndorff-Nielsen and Thorbjørnsen ([2004](#)), but the proof we give here is shorter]:

THEOREM 1.6. *Let $\mu$ be a $*$-infinitely divisible distribution with compactly supported Lévy measure, and let, for $n$ integer, $\mu_n$ be a probability measure such that $\mu_n^{*n} = \mu$. Then for each $k \geq 1$, the sequence $(n \times m_k(\mu_n))_n$ tends to $\mathfrak{C}_k(\mu)$, which is equal to $\mathfrak{K}_k(\Psi(\mu))$.*

PROOF. By [(3)](#), one has

$$n \times m_k(\mu_n) = n \sum_{\pi \in \mathrm{Part}(k)} \prod_{V \in \pi} \underbrace{\mathfrak{C}_{|V|}(\mu_n)}_{\mathfrak{C}_{|V|}(\mu)/n} = \sum_{\pi \in \mathrm{Part}(k)} n^{1-|\pi|} \mathfrak{C}_\pi(\mu) = \mathfrak{C}_k(\mu) + o(1).$$

Let us denote $\nu_n = \mu_n^{\boxplus n}$. By part 2.(f) of Theorem [1.1](#), the sequence $(\nu_n)$ converges weakly to $\Psi(\mu)$. By Hölder and Minkowski inequalities in tracial noncommutative $W^*$-probability spaces, every moment of $\nu_n$ is bounded uniformly in $n$, so the cumulants of $\nu_n$ tend to the cumulants of $\Psi(\mu)$. But by [(4)](#),

$$n \times m_k(\mu_n) = n \sum_{\pi \in \mathrm{NC}(k)} \prod_{V \in \pi} \underbrace{\mathfrak{K}_{|V|}(\mu_n)}_{\mathfrak{K}_{|V|}(\nu_n)/n} = \sum_{\pi \in \mathrm{NC}(k)} n^{1-|\pi|} \mathfrak{K}_\pi(\nu_n),$$

which tends to

$$\sum_{\pi \in \mathrm{NC}(k)} \delta_1^{|\pi|} \mathfrak{K}_\pi(\Psi(\mu)) = \mathfrak{K}_k(\Psi(\mu)). \qquad \square$$

**2. Free convolution and random matrices, choice of the models.** For $\nu$ probability measure on $\mathbb{R}$, denote by $\tilde{\nu}$, the symmetrization of $\nu$, which is the probability measure defined by $\tilde{\nu}(B) = \frac{1}{2}(\nu(B) + \nu(-B))$ for all Borel set $B$.

For $M$ Hermitian matrix, we will denote by $\mu_M$ its spectral distribution, that is, the uniform distribution on its spectrum (with multiplicity).

For $M$ complex (possibly non-Hermitian) matrix, denote by $\tilde{\mu}_{|M|}$ the symmetrization of the spectral measure of $|M|$, where $|M| = \sqrt{M^*M}$ is the unique Hermitian nonnegative matrix such that $M$ can be written $M = U|M|$, with $U$ unitary.

If $M$ is a random matrix, $\mu_M$ is a random probability mesure on the real line. For $(M_d)_d$ sequence of random matrices, we will use the notion of convergence in probability for the sequence $(\mu_{M_d})$ of random probability measures.

The rest of this section may be skipped by the reader who wants to go straight to the result. We will only explain the choice of the models, that is, is of the family's $\mathbb{P}_d^\mu$ and $\mathbb{L}_d^\mu$ of distributions.

Let us now explain in detail the choice of the family of the distributions $\mathbb{P}_d^\mu$, the distributions of the random Hermitian matrices. We would not go into as much detail for the distributions $\mathbb{L}_d^\mu$, which we construct in a similar way.

The following theorem is proved in Voiculescu ([1991]) and in Pastur and Vasilchuk ([2000]) under more restrictive hypothesis, which can easily be removed using functional calculus.

THEOREM 2.1. *Let $n$ be a positive integer. Let $\mu_1, \ldots, \mu_n$ be probability measures on $\mathbb{R}$. Let, for $d \in \mathbb{N}^*$, $(M_d^{(i)})_{i=1,\ldots,n}$ be a family of independent $d \times d$ Hermitian random matrices. We suppose that for all $i = 1, \ldots, n$, the distribution of $M_d^{(i)}$ is invariant under the unitary group's action, and $\mu_{M_d^{(i)}}$ converges in probability, when $d \to \infty$, to $\mu_i$. Then the spectral distribution of $\sum_{i=1}^n M_d^{(i)}$ converges in probability, when $d \to \infty$, to $\mu_1 \boxplus \cdots \boxplus \mu_n$.*

Let us consider a sequence $(\mu_n)$ of probability measures on $\mathbb{R}$ and a sequence $(k_n)$ of integers tending to $+\infty$ such that $\mu_n^{*k_n}$ converges weakly to a probability measure $\mu$ on $\mathbb{R}$. Let, for $n \in \mathbb{N}$, $d \in \mathbb{N}^*$, $(M_{d,n}^{(i)})_{1 \leq i \leq k_n}$ be a family of independent copies of a random Hermitian $d \times d$ matrix $M_{d,n}$, whose distribution is unitarily invariant. For every $n \in \mathbb{N}$, we suppose that $\mu_{M_{d,n}}$ converges in probability, when $d \to \infty$, to $\mu_n$.

Then we know that, for every $n \in \mathbb{N}$, the spectral distribution of $\sum_{i=1}^{k_n} M_{d,n}^{(i)}$ converges in probability, when $d \to \infty$, to $\mu_n^{\boxplus k_n}$.

Let us suppose that, on the other hand, for every $d \in \mathbb{N}^*$, $\sum_{i=1}^{k_n} M_{d,n}^{(i)}$ converges in distribution, when $n \to \infty$, to a random matrix $M_d$.

We know, by Theorem 1.1, that $\mu_n^{\boxplus k_n}$ converges, when $n \to \infty$, to the image $\Psi(\mu)$ of $\mu$ by Bercovici–Pata's bijection.

A natural question is the following: does the spectral distribution of $M_d$ converge in probability, when $d \to \infty$, to $\Psi(\mu)$?

In other words, is the limit, when $d \to \infty$, of the spectral distribution of the limit, when $n \to \infty$, of $\sum_{i=1}^{k_n} M_{d,n}^{(i)}$ equal to the limit, when $n \to \infty$, of the limit, when $d \to \infty$, of the spectral distribution of $\sum_{i=1}^{k_n} M_{d,n}^{(i)}$?

The answer of this question is affirmative in our model $[M_{d,n} = U \times \operatorname{diag}(X_{n,1}, \ldots, X_{n,d})U^*$, $U$ unitary Haar-distributed, independent of the i.i.d. random variables $X_{n,1}, \ldots, X_{n,d}$ with distribution $\mu_n]$. It can be summarized in the following diagram:

$$
\begin{array}{ccc}
M_{d,n}^{(1)} + \cdots + M_{d,n}^{(k_n)} & \xrightarrow{\ n \to \infty\ } & \mathbb{P}_d^\mu \\
\big| & & \big| \\
d \text{ goes to } \infty & & d \text{ goes to } \infty \\
\downarrow & & \downarrow \\
\text{spectral law:} & \xrightarrow{\ n \to \infty\ } & \text{spectral law:} \\
\mu_n^{\boxplus k_n} & & \Psi(\mu)
\end{array}
$$

The choice of this model is supported by the three following remarks:

1. For $d$, $n \geq 1$, if $\mu_n = \frac{1}{\pi} \frac{1/n \, dx}{(1/n)^2 + x^2}$, the expectation of the spectral distribution of $\sum_{i=1}^n M_{d,n}^{(i)}$ is $\frac{1}{\pi} \frac{dx}{1 + x^2}$.

2. For any fixed $d \geq 1$, $\sum_{i=1}^{k_n} M_{d,n}^{(i)}$ converges in distribution, when $n \to \infty$, to a distribution $\mathbb{P}_d^\mu$ which depends only on $\mu = \lim_{n \to \infty} \mu_n^{*k_n}$.

3. For every pair $(\mu, \nu)$ of $*$-infinitely divisible distributions, similarly to the relation

$$
\Psi(\mu * \nu) = \Psi(\mu) \boxplus \Psi(\nu),
$$

we have, for every $d \geq 1$,

$$
\mathbb{P}_d^\mu * \mathbb{P}_d^\nu = \mathbb{P}_d^{\mu * \nu}.
$$

This property (and the fact that all formulas depend analytically on $d$, so could be extended to noninteger $d$) opens the perspective of a continuum between the classical convolution $*$ and the free convolution $\boxplus$ for infinitely divisible measures.

Let us now explain how to construct the distributions $\mathbb{L}_d^\mu$. The following theorem is easily obtained combining the results of Haagerup and Larsen (2000) and Hiai and Petz (2000), and using functional calculus.

THEOREM 2.2.   *Let $n$ be a positive integer. Let $\mu_1, \ldots, \mu_n$ be probability measures on $\mathbb{R}$. Let, for $d \geq 1$, $(M_d^{(i)})_{i=1,\ldots,n}$ be a family of random $d \times d$* *matrices with every $M_d^{(i)}$ having a distribution invariant under the left and right actions of the unitary group. We suppose that, for every $i = 1, \ldots, n$, the distribution of $M_d^{(i)}$ is invariant under the left and right unitary group's actions, and that the symmetrization $\tilde{\mu}_{|M_d^{(i)}|}$ of the spectral distribution of $|M_d^{(i)}|$ converges in probability to $\mu_i$.*

*Then the symmetrization of the spectral distribution of*

$$\left| \sum_{i=1}^n M_d^{(i)} \right|$$

*converges in probability, when $d$ tends to infinity, to $\mu_1 \boxplus \cdots \boxplus \mu_n$.*

Let us then consider, for $\mu$ symmetric $*$-infinitely divisible distribution, a sequence $(\mu_n)$ of symmetric distributions and a sequence $(k_n)$ of integers which tends to infinity such that $\mu_n^{*k_n}$ converges weakly to $\mu$. Let $d$ be a positive integer. If for all $n$, $(M_{d,n}^{(i)})_{i=1,\ldots,k_n}$ is a family of independent copies of $U \operatorname{diag}(X_1, \ldots, X_d)V$, where $U, V, X_1, \ldots, X_d$ are independent, $U$ and $V$ are unitary Haar-distributed, and $X_1, \ldots, X_d$ are distributed according to $\mu_n$, then it appears that

$$\sum_{i=1}^{k_n} M_{d,n}^{(i)}$$

converges in distribution to a distribution $\mathbb{L}_d^\mu$ which depends only on $\mu$.

We will show that if $M_d$ is distributed according to $\mathbb{L}_d^\mu$, then $\tilde{\mu}_{|M_d|}$ converges in probability to $\Psi(\mu)$.

**3. The distributions $\mathbb{P}_d^\mu$.** $\mathbb{E}$ denotes expectation. For any distribution $\mathbb{P}$ and any function $f$ on a set of matrices, $\mathbb{E}_{\mathbb{P}}(f(M))$ denotes $\int f(M) \, d\mathbb{P}(M)$. Tr denotes the trace.

THEOREM 3.1. *Let $\mu$ be an $*$-infinitely divisible distribution. Let $(\mu_n)$ be a sequence of probability measures on $\mathbb{R}$ and $(k_n)$ a sequence of integers tending to $+\infty$ such that the sequence $\mu_n^{*k_n}$ converges weakly to $\mu$. Let, for $d \geq 1$ and $n \geq 1$, $\mathbb{Q}_d^{\mu_n}$ be the distribution of $U \operatorname{diag}(X_{n,1}, \ldots, X_{n,d})U^*$, where $U$ is a Haar-distributed unitary random matrix, independent of the $\mu_n$-distributed i.i.d. random variables $X_{n,1}, \ldots, X_{n,d}$.*

*Then the sequence $((\mathbb{Q}_d^{\mu_n})^{*k_n})$ of probability measures on the space of $d \times d$ Hermitian matrices converges weakly to a distribution $\mathbb{P}_d^\mu$.*

*Moreover, Fourier transform of the distribution $\mathbb{P}_d^\mu$ on the space of $d \times d$ Hermitian matrices with the scalar product $(M, N) \mapsto \operatorname{Tr} MN$ is given by this formula: for every Hermitian matrix $A$,*

$$(5) \qquad \mathbb{E}_{\mathbb{P}_d^\mu}(\exp(i \operatorname{Tr} AM)) = \exp(\mathbb{E}(d \times \psi_\mu(\langle u, Au \rangle))),$$

*where*

- $\psi_\mu$ *is the Lévy exponent of* $\mu$, *that is, the unique continuous function* $f$ *from* $\mathbb{R}$ *into* $\mathbb{C}$ *such that* $f(0) = 0$ *and the Fourier transform of* $\mu$ *is* $\exp \circ f$,
- $\langle \cdot, \cdot \rangle$ *is the usual Hermitian product of* $\mathbb{C}^d$,
- $u = (u_1, \ldots, u_d)$ *is a uniformly distributed random vector on the unit sphere of* $\mathbb{C}^d$.

It appears clearly that, for $\mu, \nu$ ∗-infinitely divisible distributions, $\mathbb{P}_d^\mu * \mathbb{P}_d^\nu = \mathbb{P}_d^{\mu * \nu}$.

PROOF. We will show the pointwise convergence of the Fourier transform of the distribution $(\mathbb{Q}_d^{\mu_n})^{*k_n}$ on the space of $d \times d$ Hermitian matrices. Let $A$ be a $d \times d$ Hermitian matrix with spectrum $a \in \mathbb{R}^d$. Let $F_n$ (resp. $F$) be the Fourier transform of $\mu_n^{\otimes d}$ (resp. $\mu^{\otimes d}$). Then, when $n$ tends to infinity, $k_n(F_n - 1)$ converges (uniformly on every compact set of $\mathbb{R}^d$) to the Lévy exponent $\psi$ of $\mu^{\otimes d}$ (i.e., to $\psi_\mu^{\oplus d}$, where $\psi_\mu$ is the Lévy exponent of $\mu$).

We have

$$\mathbb{E}_{(\mathbb{Q}_d^{\mu_n})^{*k_n}}(\exp(i \operatorname{Tr} AM)) = (\mathbb{E}_{\mathbb{Q}_d^{\mu_n}}(\exp(i \operatorname{Tr} AM)))^{k_n}.$$

Recall $\mathbb{Q}_d^{\mu_n}$ is invariant under the unitary action, so

$$\mathbb{E}_{(\mathbb{Q}_d^{\mu_n})^{*k_n}}(\exp(i \operatorname{Tr} AM)) = (\mathbb{E}_{\mathbb{Q}_d^{\mu_n}}(\exp(i \operatorname{Tr}(\operatorname{diag}(a)M))))^{k_n}.$$

$\mathbb{Q}_d^{\mu_n}$ is the distribution of $U \operatorname{diag}(X_{n,1}, \ldots, X_{n,d})U^*$, where $U$ is a Haar distributed unitary matrix, independent of the $\mu_n$-distributed i.i.d. random variables $X_{n,1}, \ldots, X_{n,d}$. So

$$\mathbb{E}_{(\mathbb{Q}_d^{\mu_n})^{*k_n}}(\exp(i \operatorname{Tr} AM))$$

$$= (\mathbb{E}(\exp(i \operatorname{Tr}(\operatorname{diag}(a) U \operatorname{diag}(X_{n,1}, \ldots, X_{n,d})U^*))))^{k_n}$$

$$= \left( \mathbb{E}\left( \exp\left( i \sum_{k,l=1}^d a_k X_{n,l} |u_{k,l}|^2 \right) \right) \right)^{k_n}$$

$$= \left( \mathbb{E}\left( F_n\left( \left( \sum_{k=1}^d a_k |u_{k,l}|^2 \right)_{l \in [d]} \right) \right) \right)^{k_n},$$

which can be written

$$\left( 1 + \frac{1}{k_n} \mathbb{E}\left( k_n \left( F_n\left( \left( \sum_{k=1}^d a_k |u_{k,l}|^2 \right)_{l \in [d]} \right) - 1 \right) \right) \right)^{k_n}$$

(recall $[d] = \{1, \ldots, d\}$).

But $k_n(F_n - 1)$ converges uniformly on every compact set to $\psi$ when $n \to \infty$, so we have

$$\mathbb{E}_{(\mathbb{Q}_d^{\mu n})^{*k_n}}(\exp(i\operatorname{Tr} AM)) \xrightarrow{n \to \infty} \exp\left(\mathbb{E}\left(\psi\left(\left(\sum_{k=1}^{d} a_k |u_{k,l}|^2\right)_{l \in [d]}\right)\right)\right).$$

It implies that $(\mathbb{Q}_d^{\mu n})^{*k_n}$ converges in distribution to a probability measure $\mathbb{P}_d^{\mu}$ and that the Fourier transform of $\mathbb{P}_d^{\mu}$, evaluated on a $d \times d$ Hermitian matrix $A$ with spectrum $a \in \mathbb{R}^d$, is given by

$$\mathbb{E}_{\mathbb{P}_d^{\mu}}(\exp(i\operatorname{Tr} AM)) = \exp\left(\mathbb{E}\left(\psi\left(\left(\sum_{k=1}^{d} a_k |u_{k,l}|^2\right)_{l \in [d]}\right)\right)\right).$$

But $\psi = \psi_\mu^{\oplus d}$, so

$$(6) \qquad \mathbb{E}_{\mathbb{P}_d^{\mu}}(\exp(i\operatorname{Tr} AM)) = \exp(\mathbb{E}(d \times \psi_\mu(\langle Z, a \rangle))),$$

where $\langle \cdot, \cdot \rangle$ is the usual scalar product of $\mathbb{R}^d$ and $Z = (|u_1|^2, \ldots, |u_d|^2)$, with $u = (u_1, \ldots, u_d)$ a uniformly distributed random vector on the unit sphere of $\mathbb{C}^d$.

Recall that the distribution of $u$ is invariant under the unitary action, so $\mathbb{E}(d \times \psi_\mu(\langle Z, a \rangle)) = \mathbb{E}(d \times \psi_\mu(\langle u, Au \rangle))$. □

REMARK 3.2 (The Poisson case). One can already identify $\mathbb{P}_d^{\mu}$ when $\mu = \mathcal{P}(\lambda)$ is the classical Poisson distribution with parameter $\lambda$ (denoted $\pi_{\delta_1, \lambda}^*$ in Section 1.2). It is easy, using Fourier transform, to see that, in this case, $\mathbb{P}_d^{\mu}$ is the distribution of

$$\sum_{k=1}^{X(d\lambda)} u_d(k) u_d(k)^*,$$

where $(u_d(k))_{k \geq 1}$ is an independent family of uniformly distributed random vectors on the unit sphere of $\mathbb{C}^d$, independent of the $\mathcal{P}(d\lambda)$-random variable $X(d\lambda)$.

*Explicit computation of the Fourier transform of $\mathbb{P}_d^{\mu}$—the Gaussian case.* In this section we give the distribution, the moments and the Fourier transform of the random variable $Z$ appearing in (6) of the Fourier transform of $\mathbb{P}_d^{\mu}$. In the following, we will only need the moments of $Z$.

PROPOSITION 3.3. *Let $u$ be a random vector of the unit sphere of $\mathbb{C}^d$ with uniform distribution. Then the distribution of $Z = (|u_1|^2, \ldots, |u_d|^2)$ on the* *standard d-symplexe is the uniform distribution, that is, for every bounded
Borel function $f$,*

$$\mathbb{E}(f(|u_1|^2, \ldots, |u_d|^2))$$
$$= (d-1)! \int_{x_1=0}^1 \int_{x_2=0}^{1-x_1} \cdots \int_{x_{d-1}=0}^{1-\sum_{i=1}^{d-2} x_i} f\left(x_1, \ldots, x_{d-1}, 1 - \sum_{i=1}^{d-1} x_i\right) \mathrm{d}x.$$

To prove it, write $u$ as a renormalized Gaussian standard vector on $\mathbb{C}^d$,
and do an appropriate change of variables.

We deduce, by induction on $d$, the following:

PROPOSITION 3.4. *For $d \geq 1$, for $\alpha \in \mathbb{N}^d$, denoting $s = \sum_i \alpha_i$,*

$$\mathbb{E}(|u_1|^{2\alpha_1} \cdots |u_d|^{2\alpha_d}) = (d-1)! \frac{\prod_{i=1}^d (\alpha_i!)}{(s+d-1)!} \leq (s!)^s \frac{(d-1)!}{(s+d-1)!}.$$

REMARK 3.5. When $\mu = N(0,1)$, that is, $\psi_\mu(t) = -\frac{t^2}{2}$, Proposition 3.4
allows us to compute the Fourier transform. It appears then that, when
$M_d$ has distribution $\mathbb{P}_d^\mu$, $M_d$ has the distribution of $N_d + \frac{1}{\sqrt{d+1}} X.I_d$, where
$N_d \in GUE(d, \frac{1}{d+1})$ [$GUE(d, \sigma^2)$ is the Euclidean space of Hermitian $d \times d$
matrices endowed with the standard Gaussian distribution with variance $\sigma^2$]
and $X$ is a real standard Gaussian random variable independent of $N_d$.

Proposition 3.3 allows us also to compute, by induction on $d$, the Fourier
transform of the random variable $Z$.

PROPOSITION 3.6. *Let $d \geq 2$ be an integer and let $a \in \mathbb{R}^d$ be such that
the $a_k$ are pairwise distinct. Then*

$$\mathbb{E}(\exp(i\langle a, Z \rangle)) = -(d-1)! \sum_{j=1}^d \frac{e^{ia_j}}{\prod_{k=1,\ldots,\widehat{j},\ldots,d} i(a_k - a_j)}.$$

This proposition, together with the formula

$$\mathbb{E}_{\mathbb{P}_d^\mu}(\exp(i \operatorname{Tr} AM))$$
$$= \exp\left\{ i\gamma \operatorname{Tr}(A) + d \int_{u \in \mathbb{R}} \underbrace{\left[ \mathbb{E}(e^{iu\langle Z,a \rangle}) - 1 - \frac{iu \operatorname{Tr}(A)}{d(1+u^2)} \right] \frac{1+u^2}{u^2}}_{-\frac{\mathbb{E}(\langle Z,a \rangle^2)}{2} \text{ for } u=0} \mathrm{d}G(u) \right\},$$

gives us the explicit computation of the Fourier transform of $\mathbb{P}_d^\mu$.

## 4. Convergence of the $k$th moment of the mean spectral distribution to the $k$th moment of $\Psi(\mu)$ when the Lévy measure has compact support.

### 4.1. *Statement of the result, preliminaries for the proof.*

PROPOSITION 4.1. *Let $\mu$ be an $*$-infinitely divisible distribution with compactly supported Lévy measure (in the sense of the definition given at Remark* 1.2). *Then we have*

$$\forall\, k \in \mathbb{N}, \qquad \mathbb{E}_{\mathbb{P}_d^\mu}\left(\frac{1}{d}\operatorname{Tr} M^k\right) - m_k(\Psi(\mu)) = O\left(\frac{1}{d}\right).$$

*Notation and preliminaries.* Let, for $n \in \mathbb{N}^*$, $\mu_n$ be a probability measure on $\mathbb{R}$ such that $\mu_n^{*n} = \mu$. Let us consider, for $d \geq 1$ and $n \geq 1$, $(M_{d,n}^{(i)})_{1 \leq i \leq n}$ i.i.d. random matrices with distribution $\mathbb{Q}_d^{\mu_n}$. Then we know by Theorem 3.1 that, for $d \geq 1$, the sum of the $M_{d,n}^{(i)}$'s $(i = 1, \ldots, n)$ converges in distribution to $\mathbb{P}_d^\mu$ when $n \to \infty$. We know, by Theorem 1.6, that, for all $k \in \mathbb{N}^*$, the sequence $n m_k(\mu_n)$ is bounded, and so

$$\forall\, k \geq 1, \forall\, d \geq 1 \qquad \mathbb{E}_{\mathbb{P}_d^\mu}\left(\frac{1}{d}\operatorname{Tr} M^k\right) = \lim_{n \to \infty} \mathbb{E}\left(\frac{1}{d}\operatorname{Tr}\left(\left(\sum_{i=1}^n M_{d,n}^{(i)}\right)^k\right)\right).$$

Let us then fix $k \in \mathbb{N}^*$.

### 4.2. *Computation of $\mathbb{E}_{\mathbb{P}_d^\mu}(\frac{1}{d}\operatorname{Tr} M^k)$ and proof of Proposition* 4.1. Let us define, for $d, n \geq 1$,

$$a_{d,n} = \mathbb{E}\left(\frac{1}{d}\operatorname{Tr}\left(\left(\sum_{i=1}^n M_{d,n}^{(i)}\right)^k\right)\right).$$

We have

$$a_{d,n} = \frac{1}{d}\operatorname{Tr}\left(\mathbb{E}\left(\sum_{f \in [n]^k} \prod_{r=1}^k M_{d,n}^{(f(r))}\right)\right).$$

We will transform this sum by summing on the partitions.

We denote by $\operatorname{Bij}(I)$ the set of permutations of a set $I$. Consider a partition $\pi$ of $[n]$ (we have defined $[n] = \{1, \ldots, n\}$) and $k \in [n]$. We denote by $\pi(k)$ the index of the class of $k$, after having ordered the classes according to the order of their first element [e.g., $\pi(1) = 1$; $\pi(2) = 1$ if $1 \overset{\pi}{\sim} 2$ and $\pi(2) = 2$ if $1 \overset{\pi}{\nsim} 2$]. We denote, for $l, n$ nonnegative integers, by $A_n^l$, the number of one-to-one maps from $[l]$ to $[n]$, that is, $n(n-1)\cdots(n-l+1)$.

The following lemma will be used quite often in the text.

LEMMA 4.2. *Consider $k, n \in \mathbb{N}^*$. Consider $\phi : [n]^k \to \mathbb{C}$ such that*

$$\forall f \in [n]^k, \, \forall \sigma \in \mathrm{Bij}([n]) \qquad \phi(\sigma \circ f) = \phi(f).$$

*Then*

$$\sum_{f \in [n]^k} \phi(f) = \sum_{\pi \in \mathrm{Part}(k)} A_n^{|\pi|} \phi((\pi(1), \ldots, \pi(k))).$$

By this lemma, we have

$$a_{d,n} = \frac{1}{d} \mathrm{Tr} \, \mathbb{E} \left( \sum_{\pi \in \mathrm{Part}(k)} A_n^{|\pi|} \prod_{r=1}^{k} M_{d,n}^{(\pi(r))} \right)$$

$$= \underbrace{\frac{1}{d} \mathrm{Tr} \, \mathbb{E} \left( \sum_{\pi \in \mathrm{NC}(k)} A_n^{|\pi|} \prod_{r=1}^{k} M_{d,n}^{(\pi(r))} \right)}_{\text{denoted by } v_{d,n}} + \underbrace{\frac{1}{d} \mathrm{Tr} \, \mathbb{E} \left( \sum_{\substack{\pi \in \mathrm{Part}(k) \\ \pi \notin \mathrm{NC}(k)}} A_n^{|\pi|} \prod_{r=1}^{k} M_{d,n}^{(\pi(r))} \right)}_{\text{denoted by } w_{d,n}}.$$

LEMMA 4.3. *Let $\pi$ be a partition of a totally ordered finite set $I$. Then the following assertions are equivalent:*

(i) *$\pi$ is noncrossing,*

(ii) *there exists a class $V$ of $\pi$ which is an interval, and $\pi \setminus \{V\}$ is a non-crossing partition of $I \setminus V$.*

Using several times Lemma 4.3 and integrating successively with respect to the different independent random variables, we have

$$v_{d,n} = \frac{1}{d} \mathrm{Tr} \left( \sum_{\pi \in \mathrm{NC}(k)} A_n^{|\pi|} \prod_{V \in \pi} m_{|V|}(\mu_n) \cdot I_d \right)$$

$$= \sum_{\pi \in \mathrm{NC}(k)} \underbrace{\frac{A_n^{|\pi|}}{n^{|\pi|}}}_{\xrightarrow{n \to \infty} 1} \prod_{V \in \pi} n \cdot m_{|V|}(\mu_n).$$

By Theorem 1.6, for every $k \geq 1$, one has

$$\lim_{n \to \infty} n \times m_k(\mu_n) = \mathfrak{K}_k(\Psi(\mu)).$$

So for every $d$,

$$\lim_{n \to \infty} v_{d,n} = \sum_{\pi \in \mathrm{NC}(k)} \prod_{V \in \pi} \mathfrak{K}_{|V|}(\Psi(\mu)) = m_k(\Psi(\mu)).$$

To treat the term $w_{d,n}$, let us expand the trace:

$$w_{d,n} = \frac{1}{d} \sum_{\substack{\pi \in \mathrm{Part}(k) \\ \pi \notin \mathrm{NC}(k)}} \sum_{\substack{j \in [d]^k \\ j_{k+1} := j_1}} A_n^{|\pi|} \mathbb{E}\left( \prod_{r=1}^{k} (M_{d,n}^{(\pi(r))})_{j_r, j_{r+1}} \right)$$

$$= \frac{1}{d} \sum_{\substack{\pi \in \mathrm{Part}(k) \\ \pi \notin \mathrm{NC}(k)}} A_n^{|\pi|} \sum_{\tau \in \mathrm{Part}(k)} A_d^{|\tau|} \mathbb{E}\left( \prod_{r=1}^{k} (M_{d,n}^{(\pi(r))})_{\tau(r), \tau(r+1)} \right),$$

where for each $\tau \in \mathrm{Part}(k)$, $\tau(k+1) = \tau(1)$.

Using the fact that $(M_{d,n}^{(i)})_{1 \le i \le n}$ are independent copies of a matrix with distribution $\mathbb{Q}_d^{\mu_n}$, we deduce

$$w_{d,n} = \frac{1}{d} \sum_{\substack{\pi \in \mathrm{Part}(k) \\ \pi \notin \mathrm{NC}(k)}} A_n^{|\pi|} \sum_{\tau \in \mathrm{Part}(k)} A_d^{|\tau|} \prod_{V \in \pi} \mathbb{E}_{\mathbb{Q}_d^{\mu_n}}\left( \prod_{r \in V} M_{\tau(r), \tau(r+1)} \right)$$

$$= \frac{1}{d} \sum_{\substack{\pi \in \mathrm{Part}(k) \\ \pi \notin \mathrm{NC}(k)}} A_n^{|\pi|} \sum_{\tau \in \mathrm{Part}(k)} A_d^{|\tau|} \prod_{V \in \pi} \mathbb{E}\left( \sum_{l \in [d]^V} \prod_{r \in V} (u_{\tau(r), l_r} \overline{u}_{\tau(r+1), l_r} X_{n, l_r}) \right),$$

where $U \in \mathbb{U}_d$ is Haar-distributed and independent of $(X_{n,1}, \ldots, X_{n,d})$. So, applying Lemma 4.2,

$$w_{d,n} = \frac{1}{d} \sum_{\substack{\pi \in \mathrm{Part}(k) \\ \pi \notin \mathrm{NC}(k)}} A_n^{|\pi|} \sum_{\tau \in \mathrm{Part}(k)} A_d^{|\tau|}$$

$$\times \prod_{V \in \pi} \mathbb{E}\left( \sum_{\sigma \in \mathrm{Part}(V)} A_d^{|\sigma|} \prod_{r \in V} (u_{\tau(r), \sigma(r)} \overline{u}_{\tau(r+1), \sigma(r)} X_{n, \sigma(r)}) \right)$$

integrating with respect to the $X_{n,l}$'s,

$$w_{d,n} = \frac{1}{d} \sum_{\substack{\pi \in \mathrm{Part}(k) \\ \pi \notin \mathrm{NC}(k)}} A_n^{|\pi|} \sum_{\tau \in \mathrm{Part}(k)} A_d^{|\tau|} \prod_{V \in \pi} \sum_{\sigma \in \mathrm{Part}(V)} A_d^{|\sigma|}$$

$$\times \underbrace{\mathbb{E}\left( \prod_{r \in V} (u_{\tau(r), \sigma(r)} \overline{u}_{\tau(r+1), \sigma(r)}) \right)}_{\text{denoted by } \alpha_{d, \tau, \sigma}} \prod_{v \in \sigma} m_{|v|}(\mu_n)$$

$$= \frac{1}{d} \sum_{\substack{\pi \in \mathrm{Part}(k) \\ \pi \notin \mathrm{NC}(k)}} \underbrace{\frac{A_n^{|\pi|}}{n^{|\pi|}}}_{\substack{n \to \infty \\ \longrightarrow 1}} \sum_{\tau \in \mathrm{Part}(k)} A_d^{|\tau|}$$

$$\times \prod_{V \in \pi} \sum_{\sigma \in \operatorname{Part}(V)} A_d^{|\sigma|} \underbrace{n^{1-|\sigma|}}_{\substack{n \to \infty \\ \longrightarrow \delta_1^{|\sigma|}}} \alpha_{d,\tau,\sigma} \prod_{v \in \sigma} \underbrace{n m_{|v|}(\mu_n)}_{\substack{n \to \infty \\ \longrightarrow \mathfrak{C}_{|v|}(\mu) \\ \text{by Theorem } 1.6}},$$

when $n \to \infty$, for every $\pi \notin \operatorname{NC}(k)$, for every $V \in \pi$, the only remaining $\sigma \in \operatorname{Part}(V)$ is $\{V\}$.

So one has

$$\mathbb{E}_{\mathbb{P}_d^\mu}\left(\frac{1}{d} \operatorname{Tr} M^k\right) - m_k(\Psi(\mu))$$

$$= \frac{1}{d} \sum_{\substack{\pi \in \operatorname{Part}(k) \\ \pi \notin \operatorname{NC}(k)}} \sum_{\tau \in \operatorname{Part}(k)} A_d^{|\tau|} \prod_{V \in \pi} A_d^1 \mathfrak{C}_{|V|}(\mu) \mathbb{E}\left(\prod_{r \in V} u_{\tau(r)} \overline{u}_{\tau(r+1)}\right),$$

where $u = (u_1, \ldots, u_d)$ is a uniformly distributed random vector of the unit sphere of $\mathbb{C}^d$.

Using the invariance of the distribution of $u$ under the action of diagonal unitary matrices, one sees that for all $k, l \geq 0$, $i \in [d]^k$, $j \in [d]^l$, if

$$\mathbb{E}\left(\prod_{r=1}^k u_{i_r} \prod_{r=1}^l \overline{u}_{j_r}\right) \neq 0,$$

then $k = l$ and there exists a permutation $\phi$ of $[k]$ such that for all $r$, $i_r = j_{\phi(r)}$.

So the preceding formula can be written

$$\mathbb{E}_{\mathbb{P}_d^\mu}\left(\frac{1}{d} \operatorname{Tr} M^k\right) - m_k(\Psi(\mu))$$

$$= \frac{1}{d} \sum_{\substack{\pi \in \operatorname{Part}(k) \\ \pi \notin \operatorname{NC}(k)}} \sum_{\tau \in \operatorname{acc}(\pi)} A_d^{|\tau|} d^{|\pi|} \prod_{V \in \pi} \mathfrak{C}_{|V|}(\mu) \mathbb{E}\left(\prod_{r \in V} u_{\tau(r)} \overline{u}_{\tau(r+1)}\right),$$

where for any finite totally ordered set $I$ (in which the following element of any element $x < \max I$ is denoted by $x+1$ and $\max I + 1 = \min I$), for any partition $\pi$ of $I$, $\operatorname{acc}(\pi)$ is defined to be the set of $\pi$-*acceptable* partitions, which is the set of partitions $\tau$ of $I$ such that

$$\forall V \in \pi, \exists \phi \in \operatorname{Bij}(V), \forall r \in V \qquad \tau(r) = \tau(\phi(r) + 1).$$

LEMMA 4.4.  *Let $I$ be a finite totally ordered set, $\pi, \tau$ be partitions of $I$ such that*:

- $\pi$ *has a crossing (i.e., $\pi$ is not noncrossing),*
- $\tau$ *is $\pi$-acceptable.*

*Then we have*

$$|\pi| + |\tau| \leq |I|. \tag{7}$$

PROOF. We prove the lemma by induction on the cardinality of $I$ (which is not less than four because $\pi$ has a crossing).

• If the cardinality of $I$ is four, then we can suppose $I = [4]$. We have $\pi = \{\{1,3\},\{2,4\}\}$ and the inequality (7) is easy to verify because there are only three $\pi$-acceptable partitions:

$$\{[4]\}, \qquad \{\{1,2\},\{3,4\}\}, \qquad \{\{1,4\},\{2,3\}\}.$$

• Suppose the inequality (7) proved when the cardinality of $I$ is $p$, and consider $I$ with cardinality $p+1$, and $\pi, \tau$ partitions of $I$ such that $\pi$ has a crossing and $\tau$ is $\pi$-acceptable.

• If $\pi$ and $\tau$ have no singleton class, then their cardinalities are not greater than $|I|/2$ and (7) is verified.

• If $\pi$ has a singleton class $\{a\}$, then $\tau(a) = \tau(a+1)$. This implies that if one removes the element $a$ in $I$, the class $\{a\}$ in $\pi$ and the element $a$ of its class in $\tau$, then $\tau$ stays $\pi$-acceptable (and, clearly, $\pi$ keeps a crossing). So, by induction hypothesis, we have $(|\pi| - 1) + |\tau| \leq |I| - 1$.

• If $\tau$ has a singleton class $\{b\}$, denote by $V$ the class of $b$ in $\pi$ and by $\phi$ the permutation of $V$ such that for all $r \in V$, $\tau(r) = \tau(\phi(r) + 1)$. We must have $\phi(b) + 1 = b$, so $b - 1 \overset{\pi}{\sim} b$. Remove the element $b$ in $I$, the class $\{b\}$ in $\tau$ and the element $b$ of $V$. Then, clearly, $\pi$ keeps a crossing. Define $\tilde{\phi}$ to be the permutation of the "new" $V$ by

$$\tilde{\phi}(r) = \begin{cases} \phi(r), & \text{if } \phi(r) \neq b, \\ b - 1, & \text{if } \phi(r) = b. \end{cases}$$

Then for all $r$ in the "new" $V$, $r$ and $\tilde{\phi}(r)$ are in the same class of the "new" $\tau$. It implies that $\tau$ stays $\pi$-acceptable. So, by the induction hypothesis, we have $|\pi| + (|\tau| - 1) \leq |I| - 1$. $\square$

Now recall Proposition 3.4: for $\alpha \in \mathbb{N}^d$, denoting $s = \sum_i \alpha_i$,

$$\mathbb{E}(|u_1|^{2\alpha_1} \cdots |u_d|^{2\alpha_d}) = (d-1)! \frac{\prod_{i=1}^{d}(\alpha_i!)}{(s+d-1)!} \leq (s!)^s \frac{(d-1)!}{(s+d-1)!}.$$

But for $\pi, \tau \in \mathrm{Part}(k)$, with $\tau$ $\pi$-acceptable, for all $V \in \pi$, there exists $\alpha \in \mathbb{N}^d$ such that $\sum_i \alpha_i = |V|$ and

$$\mathbb{E}\left(\prod_{r \in V} u_{\tau(r)} \overline{u}_{\tau(r+1)}\right) = \mathbb{E}(|u_1|^{2\alpha_1} \cdots |u_d|^{2\alpha_d}).$$

So, by Proposition 3.4 we have

$$
\left| \mathbb{E}_{\mathbb{P}_d^\mu} \left( \frac{1}{d} \operatorname{Tr} M^k \right) - m_k(\Psi(\mu)) \right|
$$

$$
\leq \frac{1}{d} \sum_{\substack{\pi \in \operatorname{Part}(k) \\ \pi \notin \operatorname{NC}(k)}} \sum_{\tau \in \operatorname{acc}(\pi)} A_d^{|\tau|} d^{|\pi|} |\mathfrak{C}_\pi(\mu)| \prod_{V \in \pi} (|V|!)^{|V|} \frac{(d-1)!}{(|V|+d-1)!}.
$$

Let $C$ be real such that

$$
\forall\, d \geq 1,\ \forall\, s \in [k] \qquad (s!)^s \frac{(d-1)!}{(s+d-1)!} \leq C d^{-s}.
$$

We then have

$$
\left| \mathbb{E}_{\mathbb{P}_d^\mu} \left( \frac{1}{d} \operatorname{Tr} M^k \right) - m_k(\Psi(\mu)) \right|
$$

$$
\leq \frac{1}{d} \sum_{\substack{\pi \in \operatorname{Part}(k) \\ \pi \notin \operatorname{NC}(k)}} \sum_{\tau \in \operatorname{acc}(\pi)} d^{|\tau|} d^{|\pi|} |\mathfrak{C}_\pi(\mu)| C^{|\pi|} d^{-k}.
$$

But according to (7), for all $\pi \notin \operatorname{NC}(k)$ and for all $\tau \in \operatorname{acc}(\pi)$, we have $|\tau| + |\pi| - k \leq 0$, so Proposition 4.1 is shown.

## 5. Convergence in probability of the spectral distribution to $\Psi(\mu)$ when the Lévy measure has compact support.

### 5.1. *Statement of the result and preliminaries to the proof.*

PROPOSITION 5.1. *Let $\mu$ be an $*$-infinitely divisible distribution with compactly supported Lévy measure (in the sense of the definition given at Remark 1.2). Then the spectral distribution of a random matrix with distribution $\mathbb{P}_d^\mu$ converges in probability to $\Psi(\mu)$ as $d$ tends to infinity.*

*Notation and preliminaries.* We keep the notation and the objects introduced in Section 4.1. We consider a sequence $(M_d)$ of random matrices defined on the same probability space such that for all $d$, $M_d$ has distribution $\mathbb{P}_d^\mu$, and we will prove the almost sure weak convergence of the spectral distribution of $M_d$ to $\Psi(\mu)$. It implies Proposition 5.1. Since $\Psi(\mu)$ is determined by its moments, the weak convergence of any sequence of distributions to $\Psi(\mu)$ is implied by the convergence of all moments to those of $\Psi(\mu)$.

*Let us fix* $k \geq 1$. We will show that almost surely,

$$\frac{1}{d} \operatorname{Tr} M_d^k \xrightarrow{d \to \infty} m_k(\Psi(\mu)).$$

Var denotes the variance.

Recall that by Borel–Cantelli's lemma, a sequence $(Y_d)_{d \in \mathbb{N}}$ of square-integrable real random variables converges almost surely to a real $l$ if $\sum_d (\mathbb{E}(Y_d) - l)^2$ and $\sum_d \operatorname{Var}(Y_d)$ are finite.

But we know that $\mathbb{E}_{\mathbb{P}_d^\mu}(\frac{1}{d} \operatorname{Tr} M^k) - m_k(\Psi(\mu)) = O(\frac{1}{d})$. So it suffices to show that

$$\sum_d \operatorname{Var}_{\mathbb{P}_d^\mu}\left(\frac{1}{d} \operatorname{Tr} M^k\right) < \infty.$$

We will show that $\operatorname{Var}_{\mathbb{P}_d^\mu}(\frac{1}{d} \operatorname{Tr} M^k) = O(\frac{1}{d^2})$ using the formula

$$(8) \qquad \operatorname{Var}_{\mathbb{P}_d^\mu}\left(\frac{1}{d} \operatorname{Tr} M^k\right) = \lim_{n \to \infty} \underbrace{\operatorname{Var}\left(\frac{1}{d} \operatorname{Tr}\left(\left(\sum_{i=1}^n M_{d,n}^{(i)}\right)^k\right)\right)}_{\text{denoted by } V_{d,n}}.$$

5.2. *Computation of* $\operatorname{Var}_{\mathbb{P}_d^\mu}(\frac{1}{d} \operatorname{Tr} M^k)$ *and proof of Proposition* 5.1. We have

$$V_{n,d} = \mathbb{E}\left(\sum_{f \in [n]^{2k}} \frac{1}{d} \operatorname{Tr}\left(\prod_{r=1}^k M_{d,n}^{(f(r))}\right) \frac{1}{d} \operatorname{Tr}\left(\prod_{r=k+1}^{2k} M_{d,n}^{(f(r))}\right)\right)$$

$$- \left(\mathbb{E}\left(\sum_{f \in [n]^k} \frac{1}{d} \operatorname{Tr} \prod_{r=1}^k M_{d,n}^{(f(r))}\right)\right)^2.$$

Let us apply Lemma 4.2:

$$V_{n,d} = \sum_{\pi \in \operatorname{Part}(2k)} A_n^{|\pi|} \mathbb{E}\left(\left(\frac{1}{d} \operatorname{Tr} \prod_{r=1}^k M_{d,n}^{(\pi(r))}\right)\left(\frac{1}{d} \operatorname{Tr} \prod_{r=k+1}^{2k} M_{d,n}^{(\pi(r))}\right)\right)$$

$$- \sum_{\pi_1, \pi_2 \in \operatorname{Part}(k)} A_n^{|\pi_1|} A_n^{|\pi_2|} \mathbb{E}\left(\frac{1}{d} \operatorname{Tr} \prod_{r=1}^k M_{d,n}^{(\pi_1(r))}\right) \mathbb{E}\left(\frac{1}{d} \operatorname{Tr} \prod_{r=1}^k M_{d,n}^{(\pi_2(r))}\right).$$

We split the sum into two parts: in the first one we sum over the partitions of $[2k]$ which can be split into two partitions $\pi_1$ and $\pi_2$ respectively of $\{1, \ldots, k\}$ and of $\{k+1, \ldots, 2k\}$, in the second one we sum over other partitions of $[2k]$,

$$V_{n,d} = \sum_{\pi_1, \pi_2 \in \operatorname{Part}(k)} [A_n^{|\pi_1| + |\pi_2|} - A_n^{|\pi_1|} A_n^{|\pi_2|}]$$

$$\times \mathbb{E}\left[\frac{1}{d}\operatorname{Tr}\prod_{r=1}^{k}M_{d,n}^{(\pi_1(r))}\right]\mathbb{E}\left[\frac{1}{d}\operatorname{Tr}\prod_{r=1}^{k}M_{d,n}^{(\pi_2(r))}\right]$$

$$+\sum_{\substack{\pi\in\operatorname{Part}(2k)\\\exists\, i\le k<j, i\overset{\pi}{\sim}j}}A_n^{|\pi|}\mathbb{E}\left[\frac{1}{d}\operatorname{Tr}\left(\prod_{r=1}^{k}M_{d,n}^{(\pi(r))}\right)\frac{1}{d}\operatorname{Tr}\left(\prod_{r=k+1}^{2k}M_{d,n}^{(\pi(r))}\right)\right].$$

Let us expand the trace:

$$V_{n,d}=\frac{1}{d^2}\sum_{\pi_1,\pi_2\in\operatorname{Part}(k)}(A_n^{|\pi_1|+|\pi_2|}-A_n^{|\pi_1|}A_n^{|\pi_2|})$$

$$\times\mathbb{E}\left(\sum_{\substack{j\in[d]^k\\j_{k+1}:=j_1}}\prod_{r=1}^{k}(M_{d,n}^{(\pi_1(r))})_{j_r,j_{r+1}}\right)\mathbb{E}\left(\sum_{\substack{j\in[d]^k\\j_{k+1}:=j_1}}\prod_{r=1}^{k}(M_{d,n}^{(\pi_2(r))})_{j_r,j_{r+1}}\right)$$

$$+\frac{1}{d^2}\sum_{\substack{\pi\in\operatorname{Part}(2k)\\\exists\, i\le k<j, i\overset{\pi}{\sim}j}}A_n^{|\pi|}\sum_{j\in[d]^{2k}}\mathbb{E}[(M_{d,n}^{(\pi(1))})_{j_1,j_2}\cdots(M_{d,n}^{(\pi(k))})_{j_k,j_1}$$

$$\times(M_{d,n}^{(\pi(k+1))})_{j_{k+1},j_{k+2}}\cdots(M_{d,n}^{(\pi(2k))})_{j_{2k},j_{k+1}}].$$

We apply Lemma 4.2 once more:

$$V_{n,d}=\frac{1}{d^2}\sum_{\pi_1,\pi_2\in\operatorname{Part}(k)}(A_n^{|\pi_1|+|\pi_2|}-A_n^{|\pi_1|}A_n^{|\pi_2|})\sum_{\tau_1,\tau_2\in\operatorname{Part}(k)}A_d^{|\tau_1|}A_d^{|\tau_2|}$$

$$\times\mathbb{E}\left[\prod_{r=1}^{k}(M_{d,n}^{(\pi_1(r))})_{\tau_1(r),\tau_1(r+1)}\right]\mathbb{E}\left[\prod_{r=1}^{k}(M_{d,n}^{(\pi_2(r))})_{\tau_2(r),\tau_2(r+1)}\right]$$

$$+\frac{1}{d^2}\sum_{\substack{\pi\in\operatorname{Part}(2k)\\\exists\, i\le k<j, i\overset{\pi}{\sim}j}}A_n^{|\pi|}\sum_{\tau\in\operatorname{Part}(2k)}A_d^{|\tau|}\mathbb{E}\left[\prod_{r=1}^{2k}(M_{d,n}^{(\pi(r))})_{\tau(r),\breve{\tau}(r+1)}\right],$$

where for any partition $\tau$ of $[k]$, $\tau(k+1)$ denotes $\tau(1)$, and for any partition $\tau$ of $[2k]$, $1\le r\le 2k+1$, we define

$$\breve{\tau}(r)=\begin{cases}\tau(r), & \text{if } s\notin\{k+1,2k+1\},\\\tau(1), & \text{if } r=k+1,\\\tau(k+1), & \text{if } r=2k+1.\end{cases}$$

Since $(M_{d,n}^{(i)})_{1\le i\le n}$ are independent copies of a matrix with distribution $\mathbb{Q}_d^{\mu_n}$, we have

$$V_{n,d}=\frac{1}{d^2}\sum_{\pi_1,\pi_2\in\operatorname{Part}(k)}(A_n^{|\pi_1|+|\pi_2|}-A_n^{|\pi_1|}A_n^{|\pi_2|})\sum_{\tau_1,\tau_2\in\operatorname{Part}(k)}A_d^{|\tau_1|}A_d^{|\tau_2|}$$

$$\times \prod_{V \in \pi_1} \mathbb{E}_{\mathbb{Q}_d^{\mu_n}}\left[\prod_{r \in V} M_{\tau_1(r),\tau_1(r+1)}\right] \prod_{V \in \pi_2} \mathbb{E}_{\mathbb{Q}_d^{\mu_n}}\left[\prod_{r \in V} M_{\tau_2(r),\tau_2(r+1)}\right]$$

$$+ \frac{1}{d^2} \sum_{\substack{\pi \in \mathrm{Part}(2k) \\ \exists\, i \le k < j, i \overset{\pi}{\sim} j}} A_n^{|\pi|} \sum_{\tau \in \mathrm{Part}(2k)} A_d^{|\tau|} \prod_{V \in \pi} \mathbb{E}_{\mathbb{Q}_d^{\mu_n}}\left[\prod_{r \in V} M_{\tau(r),\check\tau(r+1)}\right]$$

$$= \frac{1}{d^2} \sum_{\pi_1,\pi_2 \in \mathrm{Part}(k)} \left(A_n^{|\pi_1|+|\pi_2|} - A_n^{|\pi_1|} A_n^{|\pi_2|}\right) \sum_{\tau_1,\tau_2 \in \mathrm{Part}(k)} A_d^{|\tau_1|} A_d^{|\tau_2|}$$

$$\times \prod_{V \in \pi_1} \mathbb{E}\left[\sum_{l \in [d]^V} \prod_{r \in V} u_{\tau_1(r),l_r} \overline{u}_{\tau_1(r+1),l_r} X_{n,l_r}\right]$$

$$\times \prod_{V \in \pi_2} \mathbb{E}\left[\sum_{l \in [d]^V} \prod_{r \in V} u_{\tau_2(r),l_r} \overline{u}_{\tau_2(r+1),l_r} X_{n,l_r}\right]$$

$$+ \frac{1}{d^2} \sum_{\substack{\pi \in \mathrm{Part}(2k) \\ \exists\, i \le k < j, i \overset{\pi}{\sim} j}} A_n^{|\pi|} \sum_{\tau \in \mathrm{Part}(2k)} A_d^{|\tau|} \prod_{V \in \pi} \mathbb{E}\left[\sum_{l \in [d]^V} \prod_{r \in V} u_{\tau(r),l_r} \overline{u}_{\check\tau(r+1),l_r} X_{n,l_r}\right],$$

where $U$ is a unitary Haar-distributed random matrix, independent of $(X_{n,1},\ldots,X_{n,d})$. So, after application of Lemma 4.2,

$$V_{n,d} = \frac{1}{d^2} \sum_{\pi_1,\pi_2 \in \mathrm{Part}(k)} \left(A_n^{|\pi_1|+|\pi_2|} - A_n^{|\pi_1|} A_n^{|\pi_2|}\right) \sum_{\tau_1,\tau_2 \in \mathrm{Part}(k)} A_d^{|\tau_1|} A_d^{|\tau_2|}$$

$$\times \prod_{V \in \pi_1} \mathbb{E}\left(\sum_{\sigma \in \mathrm{Part}(V)} A_d^{|\sigma|} \prod_{r \in V} u_{\tau_1(r),\sigma(r)} \overline{u}_{\tau_1(r+1),\sigma(r)} X_{n,\sigma(r)}\right)$$

$$\times \prod_{V \in \pi_2} \mathbb{E}\left(\sum_{\sigma \in \mathrm{Part}(V)} A_d^{|\sigma|} \prod_{r \in V} u_{\tau_2(r),\sigma(r)} \overline{u}_{\tau_2(r+1),\sigma(r)} X_{n,\sigma(r)}\right)$$

$$+ \frac{1}{d^2} \sum_{\substack{\pi \in \mathrm{Part}(2k) \\ \exists\, i \le k < j, i \overset{\pi}{\sim} j}} A_n^{|\pi|} \sum_{\tau \in \mathrm{Part}(2k)} A_d^{|\tau|}$$

$$\times \prod_{V \in \pi} \mathbb{E}\left(\sum_{\sigma \in \mathrm{Part}(V)} A_d^{|\sigma|} \prod_{r \in V} u_{\tau(r),\sigma(r)} \overline{u}_{\check\tau(r+1),\sigma(r)} X_{n,\sigma(r)}\right),$$

integrating with respect to the $X_{n,l}$'s,

$$V_{n,d} = \frac{1}{d^2} \sum_{\pi_1,\pi_2 \in \mathrm{Part}(k)} \frac{A_n^{|\pi_1|+|\pi_2|} - A_n^{|\pi_1|} A_n^{|\pi_2|}}{n^{|\pi_1|+|\pi_2|}} \sum_{\tau_1,\tau_2 \in \mathrm{Part}(k)} A_d^{|\tau_1|} A_d^{|\tau_2|}$$

$$\times \prod_{V \in \pi_1} \sum_{\sigma \in \mathrm{Part}(V)} n^{1-|\sigma|} A_d^{|\sigma|} \mathbb{E}\left[\prod_{r \in V} u_{\tau_1(r),\sigma(r)} \overline{u}_{\tau_1(r+1),\sigma(r)}\right] \prod_{v \in \sigma} n \times m_{|v|}(\mu_n)$$

$$\times \prod_{V \in \pi_2} \sum_{\sigma \in \mathrm{Part}(V)} n^{1-|\sigma|} A_d^{|\sigma|} \mathbb{E}\left[\prod_{r \in V} u_{\tau_2(r),\sigma(r)} \overline{u}_{\tau_2(r+1),\sigma(r)}\right] \prod_{v \in \sigma} n \times m_{|v|}(\mu_n)$$

$$+ \frac{1}{d^2} \sum_{\substack{\pi \in \mathrm{Part}(2k) \\ \exists\, i \le k < j, i \overset{\pi}{\sim} j}} \frac{A_n^{|\pi|}}{n^{|\pi|}} \sum_{\tau \in \mathrm{Part}(2k)} A_d^{|\tau|} \prod_{V \in \pi} \sum_{\sigma \in \mathrm{Part}(V)} n^{1-|\sigma|} A_d^{|\sigma|}$$

$$\times \mathbb{E}\left[\prod_{r \in V} u_{\tau(r),\sigma(r)} \overline{u}_{\check{\tau}(r+1),\sigma(r)}\right] \prod_{v \in \sigma} n \times m_{|v|}(\mu_n).$$

Let $n$ tend to infinity:

$$V_{n,d} = \frac{1}{d^2} \sum_{\pi_1,\pi_2 \in \mathrm{Part}(k)} \underbrace{\frac{A_n^{|\pi_1|+|\pi_2|} - A_n^{|\pi_1|} A_n^{|\pi_2|}}{n^{|\pi_1|+|\pi_2|}}}_{\overset{n\to\infty}{\longrightarrow} 0} \sum_{\tau_1,\tau_2 \in \mathrm{Part}(k)} A_d^{|\tau_1|} A_d^{|\tau_2|}$$

$$\times \prod_{V \in \pi_1} \sum_{\sigma \in \mathrm{Part}(V)} \underbrace{n^{1-|\sigma|}}_{\overset{n\to\infty}{\longrightarrow} \delta_1^{|\sigma|}} A_d^{|\sigma|}$$

$$\times \mathbb{E}\left[\prod_{r \in V} u_{\tau_1(r),\sigma(r)} \overline{u}_{\tau_1(r+1),\sigma(r)}\right] \prod_{v \in \sigma} \underbrace{n \times m_{|v|}(\mu_n)}_{\substack{\overset{n\to\infty}{\longrightarrow} \mathfrak{C}_{|v|}(\mu) \\ \text{by Theorem 1.6}}}$$

$$\times \prod_{V \in \pi_2} \sum_{\sigma \in \mathrm{Part}(V)} \underbrace{n^{1-|\sigma|}}_{\overset{n\to\infty}{\longrightarrow} \delta_1^{|\sigma|}} A_d^{|\sigma|}$$

$$\times \mathbb{E}\left[\prod_{r \in V} u_{\tau_2(r),\sigma(r)} \overline{u}_{\tau_2(r+1),\sigma(r)}\right] \prod_{v \in \sigma} \underbrace{n \times m_{|v|}(\mu_n)}_{\substack{\overset{n\to\infty}{\longrightarrow} \mathfrak{C}_{|v|}(\mu) \\ \text{by Theorem 1.6}}}$$

$$+ \frac{1}{d^2} \sum_{\substack{\pi \in \mathrm{Part}(2k) \\ \exists\, i \le k < j, i \overset{\pi}{\sim} j}} \underbrace{\frac{A_n^{|\pi|}}{n^{|\pi|}}}_{\overset{n\to\infty}{\longrightarrow} 1} \sum_{\tau \in \mathrm{Part}(2k)} A_d^{|\tau|} \prod_{V \in \pi} \sum_{\sigma \in \mathrm{Part}(V)} \underbrace{n^{1-|\sigma|}}_{\overset{n\to\infty}{\longrightarrow} \delta_1^{|\sigma|}} A_d^{|\sigma|}$$

$$\times \mathbb{E}\left[\prod_{r \in V} u_{\tau(r),\sigma(r)} \overline{u}_{\check{\tau}(r+1),\sigma(r)}\right] \prod_{v \in \sigma} \underbrace{n \times m_{|v|}(\mu_n)}_{\substack{\overset{n\to\infty}{\longrightarrow} \mathfrak{C}_{|v|}(\mu) \\ \text{by Theorem 1.6}}}$$

when $n$ tends to infinity, for every partition $\pi$ (or $\pi_1$ or $\pi_2$), for every $V \in \pi$, the only resting $\sigma \in \text{Part}(V)$ is $\{V\}$. So one has

$$\text{Var}_{\mathbb{P}_d^\mu}\left(\frac{1}{d}\text{Tr}\,M^k\right)$$

$$= \frac{1}{d^2} \sum_{\substack{\pi \in \text{Part}(2k) \\ \exists i \leq k < j, i \overset{\pi}{\sim} j}} \sum_{\tau \in \text{Part}(2k)} A_d^{|\tau|} \prod_{V \in \pi} A_d^1 \mathfrak{C}_{|V|}(\mu) \mathbb{E}\left(\prod_{r \in V} u_{\tau(r)} \overline{u}_{\check{\tau}(r+1)}\right),$$

where $u = (u_1, \ldots, u_d)$ is a uniformly distributed random vector of the unit sphere of $\mathbb{C}^d$.

But recall that by invariance of the distribution of $u$ under the action of diagonal unitary matrices, for all $k, l \geq 0$, $i \in [d]^k, j \in [d]^l$, if

$$\mathbb{E}\left(\prod_{r=1}^k u_{i_r} \prod_{r=1}^l \overline{u}_{j_r}\right) \neq 0,$$

then $k = l$ and there exists a permutation $\phi$ of $[k]$ such that for all $r$, $i_r = j_{\phi(r)}$.

So the preceding formula can be written

$$\text{Var}_{\mathbb{P}_d^\mu}\left(\frac{1}{d}\text{Tr}\,M^k\right)$$

$$= \frac{1}{d^2} \sum_{\substack{\pi \in \text{Part}(2k) \\ \exists i \leq k < j, i \overset{\pi}{\sim} j}} \sum_{\tau \in \text{adm}(\pi)} A_d^{|\tau|} d^{|\pi|} \prod_{V \in \pi} \mathfrak{C}_{|V|}(\mu) \mathbb{E}\left(\prod_{r \in V} u_{\tau(r)} \overline{u}_{\check{\tau}(r+1)}\right),$$

where $\text{adm}(\pi)$ is defined in the following way (splitting the set $[2k]$ in two disjoint sets $[k]$, $[2k] \setminus [k]$): for any pair $(I, J)$ of disjoint finite totally ordered sets, for any partition $\pi$ of $I \cup J$, $\text{adm}(\pi)$ is defined to be the set of $\pi$-*admissible* partitions, which is the set of partitions $\tau$ of $I \cup J$ such that

$$\forall V \in \pi, \exists \phi \in \text{Bij}(V), \forall r \in V \qquad \tau(r) = \tau(\phi(r) + 1),$$

where for any $x \in I$ (resp. $x \in J$), $x + 1$ denotes the element following $x$ in $I$ (resp. $J$).

LEMMA 5.2.  *Let $(I, J)$ be a pair of disjoint finite totally ordered sets, $\pi, \tau$ partitions of $I \cup J$ such that:*

- *there exists $i \in I, j \in J$, with $i \overset{\pi}{\sim} j$,*
- *$\tau$ is $\pi$-admissible.*

*Then we have*

(9) $$|\pi| + |\tau| \leq |I| + |J|.$$

This inequality can be proved by induction, the proof is analaguous to the one of (7).

Recall (Proposition 3.4) that for $\alpha \in \mathbb{N}^d$, using the notation $s = \sum_i \alpha_i$,

$$(10) \quad \mathbb{E}(|u_1|^{2\alpha_1} \cdots |u_d|^{2\alpha_d}) = (d-1)! \frac{\prod_{i=1}^d (\alpha_i!)}{(s+d-1)!} \leq (s!)^s \frac{(d-1)!}{(s+d-1)!}.$$

But for every $\pi, \tau \in \mathrm{Part}(2k)$, with $\tau$ $\pi$-admissible, for every $V \in \pi$, there exists $\alpha \in \mathbb{N}^d$ such that $\sum_i \alpha_i = |V|$ and

$$\mathbb{E}\left( \prod_{r \in V} u_{\tau(r)} \overline{u}_{\bar{\tau}(r+1)} \right) = \mathbb{E}(|u_1|^{2\alpha_1} \cdots |u_d|^{2\alpha_d}).$$

So, by (10), we have

$$\mathrm{Var}_{\mathbb{P}_d^\mu}\left( \frac{1}{d} \mathrm{Tr}\, M^k \right)$$

$$\leq \frac{1}{d^2} \sum_{\substack{\pi \in \mathrm{Part}(2k) \\ \exists\, i \leq k < j, i \overset{\pi}{\sim} j \\ \tau \in \mathrm{adm}(\pi)}} A_d^{|\tau|} d^{|\pi|} \prod_{V \in \pi} \mathfrak{C}_{|V|}(\mu)(|V|!)^{|V|} \frac{(d-1)!}{(|V|+d-1)!}.$$

Consider $C < \infty$ such that

$$\forall\, d \geq 1,\, \forall\, s \in [2k] \qquad (s!)^s \frac{(d-1)!}{(s+d-1)!} \leq C d^{-s}.$$

We have

$$|V_d| \leq \frac{1}{d^2} \sum_{\substack{\pi \in \mathrm{Part}(2k) \\ \exists\, i \leq k < j, i \overset{\pi}{\sim} j \\ \tau \in \mathrm{adm}(\pi)}} d^{|\tau|} d^{|\pi|} |\mathfrak{C}_\pi(\mu)| C^{|\pi|} d^{-2k}.$$

But according to (9), for every $\pi \in \mathrm{Part}(2k)$ such that there exists $i \leq k < j$ with $i \overset{\pi}{\sim} j$ and for every $\pi$-admissible $\tau \in \mathrm{Part}(2k)$, we have $|\tau| + |\pi| - 2k \leq 0$, so

$$\mathrm{Var}_{\mathbb{P}_d^\mu}\left( \frac{1}{d} \mathrm{Tr}\, M^k \right) = O\left( \frac{1}{d^2} \right)$$

and Proposition 5.1 is proved.

5.3. *Applications to GUE and sums of independent projections.* This section is not necessary for the rest of the text.

Proposition 5.1 contains the almost sure convergence of the spectral distribution of the matrices of $GUE(d, \frac{1}{d+1})$ to the semi-circle distribution, where $GUE(d, \sigma^2)$ is the Euclidean space of $d \times d$ Hermitian matrices with the scalar product $\mathrm{Tr}(\cdot \times \cdot)$, endowed with the standard Gaussian distribution with variance $\sigma^2$.

Indeed, let $(N_d)_{d \in \mathbb{N}^*}$ be a sequence of random matrices such that for every $d$, the distribution of $N_d$ is the one of a matrix of the $GUE(d, \frac{1}{d+1})$ [we do not do any hypothesis about the joint distribution of $(N_d)_{d \in \mathbb{N}^*}$]. Let $X$ be a real Gaussian standard random variable, independent of $(N_d)_{d \in \mathbb{N}^*}$. We have seen to Remark 3.5 that for $d \in \mathbb{N}^*$, $M_d := N_d + \frac{X}{\sqrt{d+1}} \cdot I_d$ has distribution $\mathbb{P}_d^{N(0,1)}$. We have proved that $\mu_{M_d}$ converges almost surely to the centered semicircle distribution with variance 1. So $\mu_{N_d}$, which is equal to $\delta_{-\frac{X}{\sqrt{d+1}}} * \mu_{M_d}$, converges almost surely to the centered semicircle distribution with variance 1.

Another consequence of Proposition 5.1 is the following one. Recall that for all $\lambda \geq 0$, the Marchenko–Pastur distribution with index $\lambda$ is the image, by the Bercovici–Pata bijection, of the classical Poisson distribution $\mathcal{P}(\lambda)$ with index $\lambda$.

PROPOSITION 5.3.    *Let, for all $d \geq 1$, $(u_d(k))_{k \geq 1}$ be an independent family of uniformly distributed random vectors on the unit sphere of $\mathbb{C}^d$. Then for all $\lambda \geq 0$, the spectral distribution of*

$$\sum_{k=1}^{d'} u_d(k) u_d(k)^*$$

*converges in probability to the Marchenko–Pastur distribution with index $\lambda$ when $d, d'$ tend to infinity and the ratio $d'/d$ tends to $\lambda$.*

The proof of this result, which uses tools introduced in the following section, is in the Appendix.

## 6. Convergence in probability of the spectral distribution $M_d$ to $\Psi(\mu)$ without condition on the Lévy measure.

6.1. *Convergence in probability of a sequence of random distributions to a deterministic distribution.* We will denote, for $z \in \mathbb{C}$, by $\Re z$ and $\Im z$ its real and imaginary parts. Let us define, for $\nu$ probability measure on $\mathbb{R}$,

$$f_\nu : \mathbb{C}^+ = \{z \in \mathbb{C}; \Im z > 0\} \to \mathbb{C},$$

$$z \mapsto \int_{u \in \mathbb{R}} \frac{\mathrm{d}\nu(u)}{u - z}.$$

Then $f_\nu$ is a holomorphic function on $\mathbb{C}^+$, $|f_\nu(z)| \leq \frac{1}{\Im z}$, and the map

$$\{\text{probability measures on } \mathbb{R}\}^2 \to \mathbb{R}^+,$$

$$(\mu_1, \mu_2) \mapsto \sup\{|f_{\mu_1}(z) - f_{\mu_2}(z)|; \Im z \geq 1\}$$

is a distance which defines the weak topology.

So, for $(M_d)_{d \geq 1}$ sequence of Hermitian random matrices and $\rho$ probability measure on $\mathbb{R}$, we have equivalence between:

(i) the spectral distribution of $\mu_{M_d}$ converges in probability to $\rho$,

(ii) for every $\varepsilon > 0$,

$$\mathrm{P}\left( \sup_{\Im z \geq 1} \left| \frac{1}{d} \operatorname{Tr}(\mathfrak{R}_z(M_d)) - f_\rho(z) \right| > \varepsilon \right) \xrightarrow{d \to \infty} 0,$$

where, for $M$ Hermitian matrix and $z \in \mathbb{C} \setminus \mathbb{R}$, $\mathfrak{R}_z(M) = (M - z)^{-1}$.

### 6.2. *Statement of the theorem and scheme of the proof.*

**THEOREM 6.1.** *Let $\mu$ be an $*$-infinitely divisible distribution. Let, for $d \geq 1$, $M_d$ be a random matrix with distribution $\mathbb{P}_d^\mu$.*

*Then the spectral distribution $\mu_{M_d}$ of $M_d$ converges in probability to $\Psi(\mu)$.*

Scheme of the proof:

1. Notation, approximation of $M_d$ by $M_d^t$.
2. Upper bound, for $a > 0$, of $\mathrm{P}(\mathrm{rg}(N_d^t) > da)$ uniformly in $d \geq 1$.
3. Conclusion.

### 6.3. *Proof of Theorem* 6.1.

6.3.1. *Notation, approximation of $M_d$ by $M_d^t$.* Let $\gamma, G$ be such that $\mu = \nu_*^{\gamma, G}$. Recall [equation (2)] that for $t > 0$, denoting:

1. By $G_t^0$ and $G_t$ the positive finite measures on $\mathbb{R}$,

$$G_t^0(A) = G(A \cap [-t, t]), \qquad G_t(A) = G(A \setminus [-t, t])$$

for all Borel set $A$ of $\mathbb{R}$.

2. The $a_t$ the number $-\int_{u \in \mathbb{R} \setminus [-t, t]} \frac{1}{u} \, \mathrm{d}G(u)$.

3. By $\mu_t, nu_t$ the measures $\nu_*^{\gamma + a_t, G_t^0}, \nu_*^{-a_t, G_t}$,

we have the following:

(i) $\mu = \mu_t * \nu_t$, so for every $d$, $M_d$ has the distribution of $M_d^t + N_d^t$, where $M_d^t$ and $N_d^t$ are independent random matrices with respective distributions $\mathbb{P}_d^{\mu_t}$ and $\mathbb{P}_d^{\nu_t}$,

(ii) $\nu_t$ is the weak limit, when $n \to \infty$, of

$$\left( \left( 1 - \frac{\lambda_t}{n} \right) \delta_0 + \frac{\lambda_t}{n} \rho_t \right)^{*n}$$

with

$$\lambda_t = \int_{u \in \mathbb{R} \setminus [-t,t]} \frac{1 + u^2}{u^2} \, \mathrm{d}G(u), \qquad \rho_t = \frac{1}{\lambda_t} \frac{1 + u^2}{u^2} \, \mathrm{d}G_t(u).$$

So for all $d \geq 1$, the distribution $\mathbb{P}_d^{\nu_t}$ of $N_d^t$ is the weak limit of the distribution of $\sum_{i=1}^{n} N_{d,n}^{t,(i)}$, where, for every $n \geq 1$, $(N_{d,n}^{t,(i)})_{1 \leq i \leq n}$ are independent copies of $U \operatorname{diag}(X_{n,1}, \ldots, X_{n,d}) U^*$ with:

(a) $(X_{n,1}, \ldots, X_{n,d})$ i.i.d. random variables with distribution

$$\left(1 - \frac{\lambda_t}{n}\right) \delta_0 + \frac{\lambda_t}{n} \rho_t,$$

(b) $U$ unitary Haar-distributed random matrix, independent of $(X_{n,1}, \ldots, X_{n,d})$.

6.3.2. *Upper bound, for $a > 0$, of* $\mathrm{P}(\mathrm{rg}(N_d^t) > da)$ *uniformly in $d \geq 1$.* We denote by $\mathrm{rg}(M)$ the rank of a matrix $M$.

Let $a$ be a positive real.

Rank is a lower semi-continuous function, so

$$\mathbb{E}(\mathrm{rg}(N_d^t)) \leq \lim_{n \to \infty} \mathbb{E}\left(\mathrm{rg}\left(\sum_{i=1}^{n} N_{d,n}^{t,(i)}\right)\right)$$

$$\leq \lim_{n \to \infty} \mathbb{E}\left(\sum_{i=1}^{n} \mathrm{rg}(N_{d,n}^{t,(i)})\right)$$

$$= \lim_{n \to \infty} \sum_{i=1}^{n} \mathbb{E}(\mathrm{rg}(N_{d,n}^{t,(i)}))$$

$$= \lim_{n \to \infty} n \mathbb{E}(\mathrm{rg}\, U \operatorname{diag}(X_{n,1}, \ldots, X_{n,d}) U^*)$$

$$= \lim_{n \to \infty} n \sum_{l=1}^{d} \mathrm{P}(X_{n,l} \neq 0)$$

$$= \lim_{n \to \infty} nd \, \mathrm{P}(X_{n,1} \neq 0)$$

$$= \lim_{n \to \infty} nd \frac{\lambda_t}{n}$$

$$= d\lambda_t.$$

So we have

$$\mathbb{E}(\mathrm{rg}(N_d^t)) \leq d\lambda_t.$$

We deduce, with the Chebyshev inequality, that, for every $a > 0$,

$$(11) \qquad \mathrm{P}(\mathrm{rg}(N_d^t) > da) \leq \frac{1}{da} \mathbb{E}(\mathrm{rg}(N_d^t)) \leq \frac{\lambda_t}{a}.$$

6.3.3. *Conclusion.* Let $\varepsilon, \eta$ be positive reals. Let us show that there exists an integer $d_0$ such that, for every integer $d \geq d_0$,

$$\mathrm{P}\left(\sup_{\Im z \geq 1}\left|\frac{1}{d}\operatorname{Tr}(\mathfrak{R}_z(M_d)) - f_{\Psi(\mu)}(z)\right| > \varepsilon\right) \leq \eta.$$

*Choice of $t > 0$.* When $t$ tends to $+\infty$, the real $a_t$ tends to 0 and the positive finite measure $G_t^0$ converges weakly to $G$. So, by Theorem 1.1, $\nu_{\boxplus}^{\gamma + a_t, G_t^0}$ converges weakly to $\nu_{\boxplus}^{\gamma, G}$. In other words, $\Psi(\mu_t)$ converges weakly to $\Psi(\mu)$. So there exists $T_1 > 0$ such that, for all $t \geq T_1$,

$$(12) \qquad \sup_{\Im z \geq 1}|f_{\Psi(\mu_t)}(z) - f_{\Psi(\mu)}(z)| < \frac{\varepsilon}{3}.$$

When $t$ tends to $+\infty$, the real $\lambda_t$ tends to 0, so there exists $T_2 > 0$ such that, for every $t \geq T_2$,

$$(13) \qquad \lambda_t \leq \frac{\varepsilon\eta}{12}.$$

Let $t = \max(T_1, T_2)$.

For every $d \geq 1$, we have

$$\mathrm{P}\left(\sup_{\Im z \geq 1}\left|\frac{1}{d}\operatorname{Tr}(\mathfrak{R}_z(M_d)) - f_{\Psi(\mu)}(z)\right| > \varepsilon\right)$$

$$= \mathrm{P}\left(\sup_{\Im z \geq 1}\left|\frac{1}{d}\operatorname{Tr}(\mathfrak{R}_z(M_d^t + N_d^t)) - f_{\Psi(\mu)}(z)\right| > \varepsilon\right)$$

and

$$(14) \quad \begin{aligned} &\left|\frac{1}{d}\operatorname{Tr}(\mathfrak{R}_z(M_d^t + N_d^t)) - f_{\Psi(\mu)}(z)\right| \\ &\leq \left|\frac{1}{d}\operatorname{Tr}(\mathfrak{R}_z(M_d^t + N_d^t) - \mathfrak{R}_z(M_d^t))\right| \\ &\quad + \left|\frac{1}{d}\operatorname{Tr}(\mathfrak{R}_z(M_d^t)) - f_{\Psi(\mu_t)}(z)\right| + |f_{\Psi(\mu_t)}(z) - f_{\Psi(\mu)}(z)|. \end{aligned}$$

Let us deal with the *first term* of the sum (14):

We know that for every complex $d \times d$ matrix $M$, $|\frac{1}{d}\operatorname{Tr} M| \leq \frac{\|M\|}{d}\operatorname{rg}(M)$, where $\|M\|$ is the operator norm of $M$ associated to the canonical Hermitian norm on $\mathbb{C}^d$, and $\|\mathfrak{R}_z(M_d^t + N_d^t) - \mathfrak{R}_z(M_d^t)\| \leq \|\mathfrak{R}_z(M_d^t + N_d^t)\| + \|\mathfrak{R}_z(M_d^t)\| \leq \frac{2}{\Im z} \leq 2$.

Moreover, for all pair $M, N$ of Hermitian matrices, for all $z \in \mathbb{C} \setminus \mathbb{R}$, $\mathfrak{R}_z(M + N) - \mathfrak{R}_z(M) = -\mathfrak{R}_z(M+N)N\mathfrak{R}_z(M)$. So $\operatorname{rg}(\mathfrak{R}_z(M_d^t + N_d^t) - \mathfrak{R}_z(M_d^t)) \leq \operatorname{rg}(N_d^t)$.

So

$$(15) \qquad \left|\frac{1}{d}\operatorname{Tr}(\mathfrak{R}_z(M_d^t + N_d^t) - \mathfrak{R}_z(M_d^t))\right| \leq \frac{2}{d}\operatorname{rg}(N_d^t),$$

but for all $d \geq 1$,

$$\mathrm{P}\bigg(\sup_{\Im z \geq 1} \bigg| \frac{1}{d} \operatorname{Tr}(\mathfrak{R}_z(M_d^t + N_d^t) - \mathfrak{R}_z(M_d^t)) \bigg| \geq \frac{\varepsilon}{3}\bigg) \leq \mathrm{P}\bigg(\bigg| \frac{2}{d} \operatorname{rg}(N_d^t)\bigg| \geq \frac{\varepsilon}{3}\bigg)$$

$$= \mathrm{P}\bigg(\big| \operatorname{rg}(N_d^t)\big| \geq \frac{d\varepsilon}{6}\bigg)$$

$$\leq \frac{6\lambda_t}{\varepsilon} \qquad \text{by (11)}$$

$$\leq \frac{\eta}{2} \qquad \text{by (13).}$$

By inequality (12), the *third term* of the sum (14) is $\leq \frac{\varepsilon}{3}$ as soon as $\Im z \geq 1$.

Let us now deal with the *second term* of the sum (14). The Lévy measure of $\mu_t$ (in the sense of the definition given at Remark 1.2), which is $G_t^0$, is compactly supported. By Proposition 5.1 and by the other results of Section 6.1, there exists an integer $d_0$ such that, for every $d \geq d_0$,

$$\mathrm{P}\bigg(\sup_{\Im z \geq 1} \bigg| \frac{1}{d} \operatorname{Tr}(\mathfrak{R}_z(M_d^t)) - f_{\Psi(\mu_t)}(z)\bigg| \geq \frac{\varepsilon}{3}\bigg) < \frac{\eta}{2}.$$

Then for all $d \geq d_0$, replacing the terms of the sum (14) by the upper bounds we just gave, we have

$$\mathrm{P}\bigg(\sup_{\Im z \geq 1} \bigg| \frac{1}{d} \operatorname{Tr}(\mathfrak{R}_z(M_d)) - f_{\Psi(\mu)}(z)\bigg| > \varepsilon\bigg) \leq \frac{\eta}{2} + \frac{\eta}{2} + 0.$$

So, we have

$$\lim_{d \to \infty} \mathrm{P}\bigg(\sup_{\Im z \geq 1} \bigg| \frac{1}{d} \operatorname{Tr}(\mathfrak{R}_z(M_d)) - f_{\Psi(\mu)}(z)\bigg| > \varepsilon\bigg) = 0,$$

and Theorem 6.1 is proved.

## 7. Study of the non-Hermitian model.

7.1. *The distributions* $\mathbb{L}_d^\mu$. This section is the analogue, for non-Hermitian matrices, of Section 3. The distributions $\mathbb{L}_d^\mu$ are defined by the following theorem, the proof of which is analogous to the one of Theorem 3.1 using the polar decomposition of non-Hermitian matrices and the bi-unitarily invariance of the distributions $\mathbb{K}_d^{\mu_n}$.

THEOREM 7.1. *Let $\mu$ be an $*$-infinitely divisible distribution. Let $(\mu_n)$ be a sequence of probability measures on $\mathbb{R}$ and $(k_n)$ be a sequence of integers tending to $+\infty$ such that the sequence $\mu_n^{*k_n}$ converges weakly to $\mu$. Let, for $d \geq 1$ and $n \geq 1$, $\mathbb{K}_d^{\mu_n}$ be the distribution of $U\operatorname{Diag}(X_{n,1}, \ldots, X_{n,d})V$, where* $U, V$ are independent unitary Haar-distributed $d \times d$ random matrices, independent of the $\mu_n$-distributed i.i.d. random variables $X_{n,1}, \ldots, X_{n,d}$.

Then the sequence $((\mathbb{K}_d^{\mu_n})^{*k_n})$ of probability measures on the space of $d \times d$ complex matrices converges weakly to a distribution $\mathbb{L}_d^{\mu}$.

Moreover, the Fourier transform of $\mathbb{L}_d^{\mu}$ on the Euclidean space of complex $d \times d$ matrices endowed with the scalar product $(M, N) \mapsto \Re(\operatorname{Tr} M^* N)$ is given by the following formula: for all complex $d \times d$ matrix $A$,

$$(16) \qquad \mathbb{E}_{\mathbb{L}_d^{\mu}}(\exp(i\Re(\operatorname{Tr} A^* X))) = \exp(\mathbb{E}(d \times \psi_\mu(\Re(\langle u, Av\rangle)))),$$

where:

- $\psi_\mu$ is the Lévy exponent of $\mu$,
- $\langle \cdot, \cdot \rangle$ is the canonical Hermitian product of $\mathbb{C}^d$,
- $u = (u_1, \ldots, u_d)$, $v = (v_1, \ldots, v_d)$ are independent random vectors, uniformly distributed on the unit sphere of $\mathbb{C}^d$.

REMARK 7.2.   1. Notice $\mathbb{L}_d^{\mu} * \mathbb{L}_d^{\nu} = \mathbb{L}_d^{\mu * \nu}$.

2. When $\mu = N(0, 1)$, $\mathbb{L}_d^{\mu}$ is the distribution of a matrix $[M_{i,j}]$ with $(\Re M_{i,j}, \Im M_{i,j})_{1 \le i,j \le d}$ $N(0, \frac{1}{2d})$-distributed i.i.d. random variables.

3. The same construction can be done with rectangular bi-unitarily invariant random matrices. It leads, when the dimensions of the matrices tend to infinity in a certain ratio, to probability measures which are infinitely divisible with respect to a certain convolution. The studying of this convolution has led the author to construct a new noncommutative probability theory, called the rectangular free probability theory, which allows us to understand the asymptotic behavior of rectangular random matrices, as free probability theory describes the asymptotic behavior of square random matrices. It might give rise to a publication.

7.2. *Convergence of the $k$th moment to the $k$th moment of $\Psi(\mu)$ when the Lévy measure is compactly supported.*   The purpose of this section is to show the following result:

PROPOSITION 7.3.   *Let $\mu$ be a symmetric $*$-infinitely divisible distribution with compactly supported Lévy measure. Then for all integer $k$,*

$$\mathbb{E}_{\mathbb{L}_d^{\mu}}(m_k(\tilde{\mu}_{|M|})) - m_k(\Psi(\mu)) = O\left(\frac{1}{d}\right).$$

PROOF.   First, for every complex $d \times d$ matrix $M$, for all integer $k$, $m_k(\tilde{\mu}_{|M|})$ is null if $k$ is odd and is equal to $\frac{1}{d} \operatorname{Tr}(MM^*)^{k/2}$ if $k$ is even. As $\mu$ is symmetric, $\Psi(\mu)$ is symmetric. So it suffices to show that, for all $k \in \mathbb{N}^*$,

$$\mathbb{E}_{\mathbb{L}_d^{\mu}}\left(\frac{1}{d} \operatorname{Tr}(MM^*)^k\right) - m_{2k}(\Psi(\mu)) = O\left(\frac{1}{d}\right).$$

Let, for $n \in \mathbb{N}^*$, $\mu_n$ be the probability measure such that $\mu_n^{*n} = \mu$. Consider, for $d \geq 1$ and $n \geq 1$, $(M_{d,n}^{(i)})_{1 \leq i \leq n}$ i.i.d. random matrices with distribution $\mathbb{K}_d^{\mu_n}$. Then we know by Theorem 7.1, that, for every $d \geq 1$, the sum of the $M_{d,n}^{(i)}$'s $(i = 1, \ldots, n)$ converges in distribution to $\mathbb{L}_d^\mu$ when $n$ goes to $\infty$.

We know, by Theorem 1.6, that, for all $k \in \mathbb{N}^*$, the sequence $(n \times m_k(\mu_n))_n$ is bounded, and so that, for all $k, d \in \mathbb{N}^*$,

$$(17) \qquad \mathbb{E}(m_{2k}(\tilde{\mu}_{|M|})) = \lim_{n \to \infty} \mathbb{E}\left( \frac{1}{d} \mathrm{Tr} \left( \left( \sum_{i=1}^n M_{d,n}^{(i)} \right) \left( \sum_{i=1}^n M_{d,n}^{(i)*} \right) \right)^k \right).$$

*Let us fix $k \in \mathbb{N}^*$.* We are going to use (17).

Let, for $d, n \geq 1$,

$$b_{d,n} = \mathbb{E}\left( \frac{1}{d} \mathrm{Tr} \left( \left( \sum_{i=1}^n M_{d,n}^{(i)} \right) \left( \sum_{i=1}^n M_{d,n}^{(i)*} \right) \right)^k \right).$$

We have

$$b_{d,n} = \frac{1}{d} \mathrm{Tr} \left( \mathbb{E} \left( \sum_{f \in \{1, \ldots, n\}^{2k}} \prod_{r=1}^k M_{d,n}^{(f(2r-1))} M_{d,n}^{(f(2r))*} \right) \right).$$

Let us transform this sum using partitions (Lemma 4.2). Moreover, from now on, we do not write anymore the index $d$ in $M_{d,n}^{(i)}$,

$$b_{d,n} = \frac{1}{d} \mathrm{Tr} \left( \mathbb{E} \left( \sum_{\pi \in \mathrm{Part}(2k)} A_n^{|\pi|} M_n^{(\pi(1))} M_n^{(\pi(2))*} M_n^{(\pi(3))} \cdots M_n^{(\pi(2k))*} \right) \right).$$

But

$$\mathbb{E}(\underbrace{M_n^{(1)*} M_n^{(1)} \cdots M_n^{(1)*}}_{2l+1 \text{ alterned factors}}) = \mathbb{E}(\underbrace{M_n^{(1)} M_n^{(1)*} \cdots M_n^{(1)}}_{2l+1 \text{ alterned factors}}) = 0 = m_{2l+1}(\mu_n) I_d,$$

$$\mathbb{E}(\underbrace{M_n^{(1)*} M_n^{(1)} \cdots M_n^{(1)}}_{2l \text{ alterned factors}}) = \mathbb{E}(\underbrace{M_n^{(1)} M_n^{(1)*} \cdots M_n^{(1)*}}_{2l \text{ alterned factors}}) = m_{2l}(\mu_n) I_d.$$

So, for $\pi \in \mathrm{NC}(2k)$, using many times Lemma 4.3 and integrating successively with respect to the different independent random variables, we obtain

$$\frac{1}{d} \mathrm{Tr}(\mathbb{E}(M_n^{(\pi(1))*} M_n^{(\pi(2))} M_n^{(\pi(3))*} \cdots M_n^{(\pi(2k))})) = m_\pi(\mu_n).$$

Proceeding then like in Section 4.2, we show easily that

$$\mathbb{E}_{\mathbb{L}_d^\mu}(m_{2k}(\tilde{\mu}_{|M|})) - m_{2k}(\Psi(\mu))$$

$$= \frac{1}{d} \sum_{\substack{\pi, \tau \in \mathrm{Part}(2k) \\ \pi \notin \mathrm{NC}(2k)}} A_d^{|\tau|} d^{|\pi|} \mathfrak{C}_\pi(\mu) \prod_{V \in \pi} \mathbb{E}\left( \prod_{\substack{r \in V \\ r \text{ odd}}} u_{\tau(r)} v_{\tau(r+1)} \prod_{\substack{r \in V \\ r \text{ even}}} \overline{u}_{\tau(r+1)} \overline{v}_{\tau(r)} \right),$$

with $\tau(2k+1) = \tau(1)$ and where $u = (u_1, \ldots, u_d), v = (v_1, \ldots, v_d)$ are independent uniformly distributed random vectors of the unit sphere of $\mathbb{C}^d$. But as we have already seen, by invariance of the distribution of $u$ under the action of unitary diagonal matrices, for every pair $(\pi, \tau)$ of partitions of $[2k]$, if

$$\prod_{V \in \pi} \mathbb{E}\left( \prod_{\substack{r \in V \\ r \text{ odd}}} u_{\tau(r)} v_{\tau(r+1)} \prod_{\substack{r \in V \\ r \text{ even}}} \overline{u}_{\tau(r+1)} \overline{v}_{\tau(r)} \right)$$

is nonzero, then for every class $V$ of $\pi$, there exists $\phi$, permutation of $V$, which maps odd numbers to even numbers and vice versa, such that for all $r \in V$, $\tau(r) = \tau(\phi(r) + 1)$. It implies that $\tau$ is $\pi$-acceptable. Using inequality $|\tau| + |\pi| \le 2k$ [equation (7)], the inequality on the moments of a uniform random vector on the sphere of $\mathbb{C}^d$ (Proposition 3.4), we deduce, as in Section 4.2, that

$$\frac{1}{d} \sum_{\substack{\pi, \tau \in \text{Part}(2k) \\ \pi \notin \text{NC}(2k)}} A_d^{|\tau|} d^{|\pi|} \mathfrak{C}_\pi(\mu) \prod_{V \in \pi} \mathbb{E}\left( \prod_{\substack{r \in V \\ r \text{ odd}}} u_{\tau(r)} v_{\tau(r+1)} \prod_{\substack{r \in V \\ r \text{ even}}} \overline{u}_{\tau(r+1)} \overline{v}_{\tau(r)} \right) = O\left(\frac{1}{d}\right).$$

So

$$\forall k \in \mathbb{N} \qquad \mathbb{E}_{\mathbb{L}_d^\mu}(m_k(\tilde{\mu}_{|M|})) - m_k(\Psi(\mu)) = O\left(\frac{1}{d}\right). \qquad \square$$

7.3. *Convergence in probability to $\Psi(\mu)$ when the Lévy measure has compact support.* The purpose of this section is to show the following result:

PROPOSITION 7.4. *Let $\mu$ be a symmetric $*$-infinitely divisible distribution with compactly supported Lévy measure (in the sense of the definition given at Remark 1.2). Let, for each $d$, $M_d$ be a random matrix with distribution $\mathbb{L}_d^\mu$.*

*Then the symmetrization of the spectral distribution of $|M_d|$ converges weakly to $\Psi(\mu)$ when $d$ goes to infinity.*

PROOF. We will show inequalities that would imply almost sure convergence of the symmetrization of the spectral distribution of $|M_d|$ to $\Psi(\mu)$ if the matrices $M_d$ ($d \ge 1$) were defined on the same probability space. So let us suppose that the matrices are defined on the same probability space. We keep the notation and objects introduced in Section 7.2. Since $\Psi(\mu)$ is symmetric and determined by its moments, the convergence of a sequence of symmetric distributions to $\Psi(\mu)$ is implied by the convergence of all the moments of even order to those of $\Psi(\mu)$.

*Let us fix $k \geq 1$.* We will show that almost surely,

$$\frac{1}{d} \operatorname{Tr}(M_d^* M_d)^k \xrightarrow{d \to \infty} m_{2k}(\Psi(\mu)).$$

But we know that

$$\mathbb{E}_{\mathbb{L}_d^\mu}\left(\frac{1}{d} \operatorname{Tr}(MM^*)^k\right) - m_{2k}(\Psi(\mu)) = O\left(\frac{1}{d}\right).$$

So it suffices to show that

$$\operatorname{Var}_{\mathbb{L}_d^\mu}\left(\frac{1}{d} \operatorname{Tr}(MM^*)^k\right) = O\left(\frac{1}{d^2}\right).$$

We will do it using the formula

$$\operatorname{Var}_{\mathbb{L}_d^\mu}\left(\frac{1}{d} \operatorname{Tr}(MM^*)^k\right) = \lim_{n \to \infty} \operatorname{Var}\left(\frac{1}{d} \operatorname{Tr}\left(\left(\sum_{i=1}^n M_{d,n}^{(i)}\right)\left(\sum_{i=1}^n M_{d,n}^{(i)*}\right)\right)^k\right).$$

Proceeding like in Section 5.2, we obtain

$$\operatorname{Var}\left(\frac{1}{d} \operatorname{Tr}(M_d^* M_d)^k\right)$$

$$= \frac{1}{d^2} \sum_{\substack{\pi, \tau \in \operatorname{Part}(4k) \\ \exists\, i \leq 2k < j,\, i \overset{\pi}{\sim} j}} A_d^{|\tau|} d^{|\pi|} \mathfrak{C}_\pi(\mu) \prod_{V \in \pi} \mathbb{E}\left(\prod_{\substack{r \in V \\ r \text{ odd}}} u_{\tau(r)} v_{\check{\tau}(r+1)} \prod_{\substack{r \in V \\ r \text{ even}}} \overline{v}_{\tau(r)} \overline{u}_{\check{\tau}(r+1)}\right),$$

where $u = (u_1, \ldots, u_d), v = (v_1, \ldots, v_d)$ are uniformly distributed independent random vectors of the unit sphere of $\mathbb{C}^d$.

But for all couple $(\pi, \tau)$ of partitions of $[4k]$, if

$$\prod_{V \in \pi} \mathbb{E}\left(\prod_{\substack{r \in V \\ r \text{ odd}}} u_{\tau(r)} v_{\check{\tau}(r+1)} \prod_{\substack{r \in V \\ r \text{ even}}} \overline{v}_{\tau(r)} \overline{u}_{\check{\tau}(r+1)}\right)$$

is nonzero, then for all class $V$ of $\pi$, there exists a permutation $\phi$ of $V$, which maps even numbers to odd numbers and vice versa, such that for all $r \in V$, $\tau(r) = \check{\tau}(\phi(r) + 1)$, which implies that $\tau$ is $\pi$-admissible. Using the inequality $|\tau| + |\pi| \leq 4k$ [equation (9)] and the inequality on the moments of a uniform random vector on the sphere of $\mathbb{C}^d$ (Proposition 3.4), we deduce, as in Section 5, that $\operatorname{Var}_{\mathbb{L}_d^\mu}(\frac{1}{d} \operatorname{Tr}(MM^*)^k) = O(\frac{1}{d^2})$. □

REMARK 7.5. In the case where $\mu = N(0, 1)$, we have a new proof of a well-known result: the spectral distribution of a Wishart $d \times d$ matrix with $d$ degrees of freedom converges almost surely, when $d$ tends to infinity, to the distribution of $X^2$ when $X$ is a centered semi-circular random variable with variance 1, which is [see Speicher (1999)] the Marchenko–Pastur distribution with parameter 1.

7.4. *Convergence in probability of $\tilde{\mu}_{|M_d|}$ to $\Psi(\mu)$ in the general case.*

THEOREM 7.6. *Let $\mu$ be a symmetric $*$-infinitely divisible distribution. Let, for $d \geq 1$, $M_d$ be a random matrix with distribution $\mathbb{L}_d^\mu$.*

*Then the symmetrization $\tilde{\mu}_{|M_d|}$ of the spectral distribution of $|M_d|$ converges in probability to $\Psi(\mu)$.*

The proof is quite similar to the one of Theorem 6.1.

PROOF OF THEOREM 7.6.

*Notation, approximation of $M_d$ by $M_d^t$.* Let $G$ be the symmetric positive finite measure on $\mathbb{R}$ such that $\mu = \nu_*^{0,G}$. Recall [equation (2)] that, for $t > 0$, if:

1. $G_t^0$ and $G_t$ are the positive finite measures

$$G_t^0(A) = G(A \cap [-t, t]), \qquad G_t(A) = G(A \setminus [-t, t])$$

   for all Borel set $A$ of $\mathbb{R}$.
2. $\mu_t = \nu_*^{0,G_t^0}$ and $\nu_t = \nu_*^{0,G_t}$,

then we have:

(i)  $\mu = \mu_t * \nu_t$, so for all $d \geq 1$, $M_d$ has the same distribution as $M_d^t + N_d^t$, where $M_d^t$ and $N_d^t$ are independent random matrices with respective distributions $\mathbb{L}_d^{\mu_t}$ et $\mathbb{L}_d^{\nu_t}$,

(ii)  $\nu_t$ is the weak limit, when $n \to \infty$, of

$$\left( \left(1 - \frac{\lambda_t}{n}\right)\delta_0 + \frac{\lambda_t}{n}\rho_t \right)^{*n}$$

for

$$\lambda_t = \int_{u \in \mathbb{R} \setminus [-t,t]} \frac{1 + u^2}{u^2} \, \mathrm{d}G(u), \qquad \rho_t = \frac{1}{\lambda_t} \frac{1 + u^2}{u^2} \, \mathrm{d}G_t(u).$$

So for all $d \geq 1$, the distribution $\mathbb{L}_d^{\nu_t}$ of $N_d^t$ is the weak limit of the distribution of $\sum_{i=1}^n N_{d,n}^{t,(i)}$, where for all $n \geq 1$, $(N_{d,n}^{t,(i)})_{1 \leq i \leq n}$ are independent copies of $U \operatorname{diag}(X_{n,1}, \ldots, X_{n,d}) V$ with:

(a)  $(X_{n,1}, \ldots, X_{n,d})$ i.i.d. random variables with distribution $(1 - \frac{\lambda_t}{n})\delta_0 + \frac{\lambda_t}{n}\rho_t$,

(b)  $U, V$ unitary Haar-distributed random matrices, independent of $(X_{n,1}, \ldots, X_{n,d})$.

In the same way as in Section 6.3.2, we show that, for all $a > 0$,

$$(18) \qquad \mathrm{P}(\mathrm{rg}(N_d^t) > da) \le \frac{\lambda_t}{a}.$$

Let us denote, for $\rho$ probability measure on $\mathbb{R}$, $\rho^2$ the distribution of $X^2$ when $X$ is a random variable with distribution $\rho$.

Consider $\varepsilon, \eta > 0$.

Let us show that there exists an integer $d_0$ such that, for all $d \ge d_0$,

$$\mathrm{P}\left(\sup_{\Im z \ge 1} \left| \frac{1}{d} \mathrm{Tr}(\mathfrak{R}_z(M_d^* M_d)) - f_{\Psi(\mu)^2}(z) \right| > \varepsilon \right) \le \eta.$$

*Choice of $t > 0$.*   When $t$ tends to $+\infty$, the measure $G_t^0$ converges weakly to $G$. So, by Theorem 1.1, $\nu_\boxplus^{0,G_t^0}$ converges weakly to $\nu_\boxplus^{0,G}$. In other words, $\Psi(\mu_t)$ converges weakly to $\Psi(\mu)$. So $\Psi(\mu_t)^2$ converges weakly to $\Psi(\mu)^2$. Hence, there exists $T_1 > 0$ such that, for all $t \ge T_1$,

$$(19) \qquad \sup_{\Im z \ge 1} |f_{\Psi(\mu_t)^2}(z) - f_{\Psi(\mu)^2}(z)| < \frac{\varepsilon}{3}.$$

When $t$ tends to $+\infty$, the real $\lambda_t$ tends to 0, so there exists $T_2 > 0$ such that, for all $t \ge T_2$,

$$(20) \qquad \lambda_t \le \frac{\varepsilon \eta}{24}.$$

Let $t = \max(T_1, T_2)$.

For all $d \ge 1$, we have

$$\mathrm{P}\left(\sup_{\Im z \ge 1} \left| \frac{1}{d} \mathrm{Tr}(\mathfrak{R}_z(M_d^* M_d)) - f_{\Psi(\mu)^2}(z) \right| > \varepsilon \right)$$

$$= \mathrm{P}\left(\sup_{\Im z \ge 1} \left| \frac{1}{d} \mathrm{Tr}(\mathfrak{R}_z((M_d^{t*} + N_d^{t*})(M_d^t + N_d^t))) - f_{\Psi(\mu)}(z) \right| > \varepsilon \right).$$

Hence,

$$\left| \frac{1}{d} \mathrm{Tr}(\mathfrak{R}_z((M_d^{t*} + N_d^{t*})(M_d^t + N_d^t))) - f_{\Psi(\mu)}(z) \right|$$

$$\le \left| \frac{1}{d} \mathrm{Tr}(\mathfrak{R}_z((M_d^{t*} + N_d^{t*})(M_d^t + N_d^t)) - \mathfrak{R}_z(M_d^{t*} M_d^t)) \right|$$

$$+ \left| \frac{1}{d} \mathrm{Tr}(\mathfrak{R}_z(M_d^{t*} M_d^t)) - f_{\Psi(\mu_t)^2}(z) \right| + |f_{\Psi(\mu_t)^2}(z) - f_{\Psi(\mu)^2}(z)|.$$

But for all pair $(M, N)$ of Hermitian matrices, for all $z \in \mathbb{C} \setminus \mathbb{R}$, $\mathfrak{R}_z(M + N) - \mathfrak{R}_z(M) = -\mathfrak{R}_z(M + N)N\mathfrak{R}_z(M)$.

So denoting $\Delta_d^t = M_d^{t*}N_d^t + N_d^{t*}(M_d^t + N_d^t)$, we have

$$
\begin{aligned}
(21) \quad & \left| \frac{1}{d} \operatorname{Tr}(\mathfrak{R}_z((M_d^{t*} + N_d^{t*})(M_d^t + N_d^t))) - f_{\Psi(\mu)}(z) \right| \\
& \leq \left| \frac{1}{d} \operatorname{Tr}(\mathfrak{R}_z((M_d^{t*} + N_d^{t*})(M_d^t + N_d^t))\Delta_d^t \mathfrak{R}_z(M_d^{t*}M_d^t)) \right| \\
& \quad + \left| \frac{1}{d} \operatorname{Tr}(\mathfrak{R}_z(M_d^{t*}M_d^t)) - f_{\Psi(\mu_t)^2}(z) \right| + |f_{\Psi(\mu_t)^2}(z) - f_{\Psi(\mu)^2}(z)|.
\end{aligned}
$$

The conclusion is similar to the one of Section 6.3.3. $\quad\square$

## APPENDIX

We prove Proposition 5.3. Denote $\sum_{k=1}^{d'} u_d(k)u_d(k)^* = N_{d,d'}$. As explained in Section 6.1, it is equivalent to prove that, for every $\varepsilon > 0$,

$$
\mathrm{P}\left( \sup_{\Im z \geq 1} \left| \frac{1}{d} \operatorname{Tr}(\mathfrak{R}_z(N_{d,d'})) - f_{\Psi(\mathcal{P}(\lambda))}(z) \right| > \varepsilon \right) \xrightarrow[d'/d \simeq \lambda]{d,d' \to \infty} 0.
$$

By Proposition 5.1, $\Psi(\mathcal{P}(\lambda))$ is the limit of the spectral distribution of a random matrix with distribution $\mathbb{P}_d^{\mathcal{P}(\lambda)}$. But, as noticed in Remark 3.2, $\mathbb{P}_d^{\mathcal{P}(\lambda)}$ is the distribution of

$$
M_d := \sum_{k=1}^{X(d\lambda)} u_d(k)u_d(k)^*,
$$

where $X(d\lambda)$ is a $\mathcal{P}(d\lambda)$-random variable, independent of the sequence $(u_d(k))$. So it suffices to prove, that for every $\varepsilon > 0$,

$$
\mathrm{P}\left( \sup_{\Im z \geq 1} \left| \frac{1}{d} \operatorname{Tr}(\mathfrak{R}_z(N_{d,d'}) - \mathfrak{R}_z(M_d)) \right| > \varepsilon \right) \xrightarrow[d'/d \simeq \lambda]{d,d' \to \infty} 0.
$$

But it was noticed in Section 6.3.3, equation (15), that

$$
\left| \frac{1}{d} \operatorname{Tr}(\mathfrak{R}_z(N_{d,d'}) - \mathfrak{R}_z(M_d)) \right| \leq \frac{2}{d} \operatorname{rg}(N_{d,d'} - M_d),
$$

which is not greater than $\frac{2}{d}|X(d\lambda) - d'|$, which converges in probability to zero, by the weak law of large numbers.

**Aknowledgments.** The author would like to thank the referees for contructive comments, and his advisor Philippe Biane, as well as Leonid Pastur, Thierry Cabanal-Duvillard and Michael Anshelevich for many interesting discussions. Also, he would like to thank Cécile Martineau for her contribution to the English version of this paper.

## REFERENCES

ANSHELEVICH, M. V. (2001). Partition-dependent stochastic measures and $q$-deformed cumulants. *Doc. Math.* **6** 343–384. MR1871667

BARNDORFF-NIELSEN, O. E. and THORBJØRNSEN, S. (2002). Selfdecomposability and Lévy processes in free probability. *Bernoulli* **3** 323–366. MR1913111

BARNDORFF-NIELSEN, O. E. and THORBJØRNSEN, S. (2004). A connection between free and classical infinite divisibility. *Inf. Dimens. Anal. Quantum Probab. Relat. Top.* **7** 573–590. MR2105912

BERCOVICI, H., PATA, V., with an appendix by BIANE, P. (1999). Stable laws and domains of attraction in free probability theory. *Ann. of Math.* **149** 1023–1060. MR1709310

BERCOVICI, H. and VOICULESCU, D. (1993). Free convolution of measures with unbounded supports. *Indiana Univ. Math. J.* **42** 733–773. MR1254116

CABANAL-DUVILLARD, T. (2004). About a matricial representation of the Bercovici–Pata bijection: A Lévy processes approach. Preprint. Available at www.math-info.univ-paris5.fr/~cabanal/liste-publi.html.

GNEDENKO, V. and KOLMOGOROV, A. N. (1954). *Limit Distributions for Sums of Independent Random Variables.* Adisson–Wesley, Reading, MA.

HAAGERUP, U. and LARSEN, F. (2000). Brown's spectral distribution measure for R-diagonal elements in finite von Neumann algebras. *J. Funct. Anal.* **176** 331–367. MR1784419

HIAI, F. and PETZ, D. (2000). *The Semicircle Law, Free Random Variables, and Entropy.* Amer. Math. Soc., Providence, RI. MR1746976

PASTUR, L. and VASILCHUK, V. (2000). On the law of addition of random matrices. *Comm. Math. Phys.* **214** 249–286. MR1796022

PETROV, V. V. (1995). *Limit Theorems of Probability Theory.* Clarendon Press, Oxford. MR1353441

SPEICHER, R. (1994). Multiplicative functions on the lattice of non-crossing partitions and free convolution. *Math. Ann.* **298** 611–628. MR1268597

SPEICHER, R. (1999). Notes of my lectures on combinatorics of free probability. Available at www.mast.queensu.ca/~speicher.

VOICULESCU, D. V. (1991). Limit laws for random matrices and free products. *Invent. Math.* **104** 201–220. MR1094052

DÉPARTEMENT DE MATHÉMATIQUES
ET APPLICATIONS
ÉCOLE NORMALE SUPÉRIEURE
45 RUE D'ULM
75230 PARIS CEDEX 05
FRANCE
E-MAIL: benaych@dma.ens.fr
URL: www.dma.ens.fr/~benaych